\documentclass[preprint,11pt]{elsarticle}

\usepackage{psfrag}
\usepackage{color}
\usepackage{amsmath}
\usepackage{amsfonts}
\usepackage{graphicx}
\usepackage{amssymb}
\usepackage{rotating}
\usepackage{multirow}
\usepackage{todonotes}
\usepackage{relsize}
\usepackage[utf8]{inputenc}
\usepackage{siunitx}
\usepackage{xargs}

\graphicspath{{FIGURES/}}

\usepackage{geometry}
\def\qed{\hbox{${\vcenter{\vbox{                        %HOLLOW SQUARE
   \hrule height 0.4pt\hbox{\vrule width 0.4pt height 6pt
   \kern5pt\vrule width 0.4pt}\hrule height 0.4pt}}}$}}
\newtheorem{prop}{Proposition}[section]
\newtheorem{rem}{Remark}[section]

\newtheorem{coroll}{Corollary}[section]

\renewcommand\v{\boldsymbol{v}}
\renewcommand\u{\boldsymbol{u}}
\newcommand\f{\boldsymbol{f}}

\newcommand\w{\boldsymbol{w}}

\renewcommand\div {{\rm div}\,}

\def\M{{\mathcal M}}

\def\v{{\bf v}}
\def\bV{{\bf V}}
\def\u{{\bf u}}
\def\f{{\bf f}}

\def\n{{\bf n}}

\def\RR{{\mathbb R}}

\def\a{{\bf a}}
\def\b{{\bf b}}

\def\R{\mathbb{R}}

\definecolor{blue}{rgb}{0,0,0.99}

\definecolor{dred}{rgb}{0.92,0,0}

\newcommand\x{\boldsymbol{x}}
\newcommand\p{\boldsymbol{p}}
\newcommand\q{\boldsymbol{q}}
\renewcommand\a{\boldsymbol{a}}
\renewcommand\b{\boldsymbol{b}}

\newcommand\m{\mathbf{m}}
\newcommand{\mId}{\mathbf{m_I}}

\def\D{{\mathcal D}}
\def\G{{\mathcal G}}
\def\Poly{{\mathbb P}}
\newcommand{\VFP}{\boldsymbol{V}_{k,h}^{\text{face}}(P)}
\newcommand{\VPP}{Q_{k-1,h}(P)}
\newcommand{\VFG}{\boldsymbol{V}_{k,h}(\Omega_h)}
\newcommand{\VPG}{Q_{k-1,h}(\Omega_h)}

\newcommandx{\GPerp}[1][1=k, usedefault]{\G_{#1}^\oplus(P)}

\newcommand{\aalpha}{\boldsymbol{\alpha}}
\newcommand{\bbeta}{\boldsymbol{\beta}}
\newcommand{\ggamma}{\boldsymbol{\gamma}}
\newcommand{\ddelta}{\boldsymbol{\delta}}
\newcommand{\bPi}{\boldsymbol{\Pi}}

\newcommand{\bcurl}{\textbf{curl}}

\newcommand{\df}{\,\text{d}f}
\newcommand{\dP}{\,\text{d}P}

\newcommand{\Cube}{\textbf{Cube}}
\newcommand{\Octa}{\textbf{Octa}}
\newcommand{\CVT}{\textbf{CVT}}
\newcommand{\Random}{\textbf{Random}}

\newcommand{\proof}{\noindent\textbf{Proof.}\hspace{1em}}
\newcommand{\cvd}{\begin{flushright}
$\square$                   
\end{flushright}}

\begin{document}

\begin{frontmatter}

\title{Parallel solvers for virtual element discretizations \\ of elliptic equations in mixed form}

\author[franco]{F. Dassi}
\ead{franco.dassi@unimib.it}
\author[simone]{S. Scacchi}
\ead{simone.scacchi@unimi.it}
\address[franco]{Dipartimento di Matematica, Universit\`a degli Studi di Milano-Bicocca, Via Cozzi, 20133 Milano, Italy}
\address[simone]{Dipartimento di Matematica, Universit\`a degli Studi di Milano, Via Saldini 50, 20133 Milano, Italy}

\begin{abstract}
The aim of this paper is twofold.
On the one hand, we test numerically the performance of mixed virtual elements in three dimensions
for the first time in the literature 
to solve the mixed formulation of three-dimensional elliptic equations on polyhedral meshes.
On the other hand, we focus on the parallel solution of the linear system arising from such discretization, 
considering both direct and iterative parallel solvers.
In the latter case, we develop two block preconditioners, one based on the approximate
Schur complement and one on a regularization technique.
Both these topics are numerically validated by several parallel tests performed on a Linux cluster. 
More specifically, 
we show that the proposed VEM discretization recovers the expected theoretical convergence properties and 
we analize the performance of the direct and iterative parallel solvers taken into account.
\end{abstract}

\end{frontmatter}

\medskip

%-------------------------------------------------------------
\section{Introduction}

In recent years, the interest of an increasing number
of researchers has focused on the development of numerical methods 
for the approximation of partial differential equations (PDEs) on polygonal or
polyhedral grids,
see e.g. \cite{beiraoE.2016}.
Among the different methodologies, the Virtual
Element Method (VEM), introduced in the pioneering paper \cite{volley},
represents a generalization of the Finite Element Method that can easily handle
general polytopal meshes. VEM can be
regarded as an evolution of the Mimetic Finite Difference method, see e.g. \cite{beiraoLM.2014}.
So far, VEM has been analyzed for general elliptic problems \cite{mixedVEM2d}, 
elasticity \cite{beiraoBM.2013,gainTP.2014}, Cahn-Hilliard \cite{ABSVu}, Stokes \cite{beiraoLV.2017},
%\cite{antoniettiBMV.2014,beiraoLV.2017},
parabolic and hyperbolic equations \cite{vaccaB.2015,vacca.2017}, discrete fracture networks \cite{benedettoBPS.2014,fumagalli2018dual} and several further applications.
Different variants of the VEM have been proposed
and analysed: 
H(div) and H(curl)-conforming \cite{divCurl}, 
%$hp$ \cite{beiraoCMR.2016}, 
serendipity \cite{chinaSere} and nonconforming \cite{ayusoLM.2016,cangianiMS.2016} VEM.

VEM for mixed formulation of two-dimensional elliptic problems has been developed
in \cite{brezziFM.2014,mixedVEM2d}; for alternative polyhedral discretizations of elliptic equations 
in mixed form see~\cite{brezziLipMixed}. We recall that second or fourth order elliptic problems 
can be reformulated as a coupled system by introducing a new variable, typically the
gradient of the potential function. This mixed formulation yields a higher order of accuracy
for the new variable and has also the favorable property of local mass conservation. 
However, similarly to other saddle point problems, the resulting linear systems are 
highly ill-conditioned due to some coefficients, such as a diffusion coefficient,
taking widely varying values in the domain where the PDE is posed. 
Consequently, the solution of saddle point matrix equations with iterative methods
(see e.g. \cite{saad.2003,zulehner.2001}) requires the construction of robust and effective preconditioners.
Successful preconditioners are based on approximate block factorization, see e.g. \cite{axelsson.1994,golubG.2003,simoncini.2004,benziGL.2005,mardal.2011,axelssonBBKA.2016}.
Several Domain Decomposition preconditioners have also been developed for finite element discretizations of
such problems, see \cite{klawonn.1998,mathew.1993a,mathew.1993b,ohWZD.2017,zampiniT.2017}.

To our knowledge, in the VEM literature only a few studies have focused on the conditioning of the stiffness matrix resulting from VEM discretizations (see \cite{mascotto.2018,dassiM.2018}) and on the development of preconditioners for VEM approximations of PDEs (see \cite{bertoluzzaPP.2017,antoniettiMasV.2018,calvo.2018}). 
Preconditioners for other polyhedral discretizations have been studied in~\cite{antoVeraZika,antonietti2015multigrid}.
We remark that all these works concern scalar elliptic equations in primal form.

The aim of the present contribution is twofold. First, we numerically verify for the first time 
that the convergence of mixed VEM scheme is in agreement with the expected theoretical estimates. Then, we develop a parallel solver for the solution of the linear systems arising from the discretization process. In the design of the parallel solver, we consider two block preconditioners, one based on the approximate Schur complement and one on a regularization technique, see e.g. \cite{benziGL.2005} and \cite{axelssonBBKA.2016}, respectively. We compare the iterative methods against the parallel direct solver Mumps \cite{amestoy.2001,amestoy.2006} by performing several parallel tests on a Linux cluster with varying number of processors, order of VEM discretization and type of polyhedral grid.

The remand of the paper is organized as follows.
In Section~\ref{pol_spaces} we give the basic notation and give some useful results on polynomial spaces and decompositions. 
Then, in Section~\ref{mod_prob}, we briefly introduce the variational formulation of the model problem.
In Section~\ref{vem_disc} we provide the Virtual Element approximation of such problem and some numerical results on the convergence of the proposed schemes.
Finally, in Section~\ref{para_precon} we focus on the parallel implementation of the discretization process, we describe the block preconditioners used for the solution of the linear system and we show some numerical experiments on a Linux cluster.

%-------------------------------------------------------------
\section{Polynomial spaces and bases}
\label{pol_spaces}

In order to develop the VEM approximation of the 3d elliptic equations in mixed form, 
we need to define suitable basis functions and a polynomial decomposition for certain
polynomial spaces introduced below.

Let $\Poly_k(\D)$ be the space of polynomials of degree $k$ defined in a general domain $\D$.
There are several choices of basis of $\Poly_k(\D)$. 
In the virtual element framework it is convenient to consider the so-called \emph{scaled monomial basis}~\cite{autostoppisti}.
Defining the multi-index $\aalpha:=(\alpha_1,\,\alpha_2,\,\alpha_3)$ with the usual notation 
$|\aalpha|:=\alpha_1+\alpha_2+\alpha_3$,
then a generic \emph{scaled monomial} on a polyhedron $P$ is 
\begin{equation*}
m_{\aalpha} := \left(\frac{\x-\x_P}{h_P}\right)^{\aalpha} = 
\left(\frac{x-x_P}{h_P}\right)^{\alpha_1}
\left(\frac{y-y_P}{h_P}\right)^{\alpha_2}
\left(\frac{z-z_P}{h_P}\right)^{\alpha_3}\,,
\end{equation*}
where $\x_P$ is the barycenter of $P$ whose coordinates are $(x_P,\,y_P,\,z_P)$ and $h_P$ is the diameter of $P$.
It is easy to show that the set of scaled monomials
\begin{equation}
\M_k(P) :=\left\{m_{\aalpha}\::\: 0\leq|\aalpha|\leq k \right\}\,, 
\label{eqn:monoBasis}
\end{equation}
is a basis for $\Poly_k(P)$.

Since we are going to exploit polynomials defined on a face $f$ of $P$,
it will be useful to introduce 2d scaled monomials. 
In such case we have to consider a multi-index $\bbeta:=(\beta_1,\,\beta_2)$ composed by only two components.
Then, a generic scaled monomial on a face $f$ is 
\begin{equation*}
m_{\bbeta}^f := \left(\frac{\widetilde{\x}-\widetilde{\x}_f}{h_f}\right)^{\bbeta} = 
\left(\frac{\widetilde{x}-\widetilde{x}_f}{h_f}\right)^{\beta_1}
\left(\frac{\widetilde{y}-\widetilde{y}_f}{h_f}\right)^{\beta_2}\,,
\end{equation*}
where $\widetilde{\x}_f=(\widetilde{x}_f,\,\widetilde{y}_f)$ 
is the barycenter of the face $f$ written in the face local coordinates system $\widetilde{x}O\widetilde{y}$ and 
$h_f$ is the diameter of the face $f$. 
As in the three dimensional case, 
the set
$$
\M_k(f) :=\left\{m_{\bbeta}^f\::\: 0\leq|\bbeta|\leq k \right\}\,, 
$$
is a basis for $\Poly_k(f)$.

Starting from such polynomial basis, 
it is possible to define also a vectorial monomial basis for the polynomial space $[\Poly_k(\D)]^d$.
We refer to this basis as $[\M_k(P)]^3$ and $[\M_k(f)]^2$ for the three and two dimensional case, respectively.

\subsection{Polynomial decomposition}

In this subsection we introduce a polynomial decomposition 
exploited to build the projection operator and the discrete forms of the model problem,
see Subsections~\ref{sub:velSpace} and~\ref{sub:discForms}.
Let us consider a polyhedron $P$, the vectorial polynomial space $[\Poly_k(P)]^3$ 
can be split in a direct sum of two spaces
\begin{equation}
[\Poly_k(P)]^3 = \G_k(P) \oplus \GPerp\,.
\label{eqn:spaceDec}
\end{equation}
where 
$$
\G_k(P) := \left\{\p_k\in[\Poly_k(P)]^3\::\: \exists\,p_{k+1}\in\Poly_{k+1}(P)\text{ such that } \p_k = \nabla p_{k+1}\right\}\,,
$$
and $\GPerp$ is complement orthogonal to $\G_k(P)$ in $[\Poly_k(P)]^3$.
A direct consequence of Equation~\eqref{eqn:spaceDec} is that a generic vectorial polynomial $\p_k\in[\Poly_k(P)]^3$ can be written as
\begin{equation}
\p_k = \nabla q_{k+1} + \x \wedge \q_{k-1}\,,
\label{eqn:decomp}
\end{equation}
where $q_{k+1}\in\Poly_{k+1}(P)$, $\q_{k-1}\in[\Poly_{k-1}(P)]^3$ and $\x:=(x,\,y,\,z)^T$ and 
we further underline that $\nabla~q_{k+1}~\in~\mathcal{G}_k(P)$ and
$\x \wedge \q_{k-1}\in\GPerp$~\cite{maxwellGO}.

Finding $q_{k+1}$ and $\q_{k-1}$ in Equation~\eqref{eqn:decomp} is not an easy task.
However, if we are dealing with vectorial scaled monomials, 
we found a straightforward recipe to get such decomposition.
From now on we will consider a multi-index $\aalpha = (\alpha_1,\,\alpha_2,\,\alpha_3)$ and 
we define the scaled vectorial monomial
$$
\mId := \left(\frac{x-x_P}{h_P},\,\frac{y-y_P}{h_P},\,\frac{z-z_P}{h_P}\right)^T\,.
$$
%----------------------------------------------
%              PRIMA COMPONENTE
%----------------------------------------------
\begin{prop}
Considering a vectorial scaled monomial with only the \emph{first} component different from 0,
Equation~\eqref{eqn:decomp} becomes
\begin{equation}
\left(
\begin{array}{c}
m_{\aalpha}\\
0 \\ 
0 
\end{array}\right) = 
\left(\frac{h_P}{|\aalpha|+1}\right) \nabla m_{\bbeta} - \left(\frac{\alpha_3}{|\aalpha|+1}\right) \mId \wedge
\left(
\begin{array}{c}
0 \\
m_{\ggamma} \\ 
0
\end{array}\right) + \left(\frac{\alpha_2}{|\aalpha|+1}\right) \mId \wedge
\left(
\begin{array}{c}
0 \\
0 \\ 
m_{\ddelta}
\end{array}\right)
\label{eqn:decForX}
\end{equation}
where 
$$
\bbeta = (\alpha_1+1,\,\alpha_2,\,\alpha_3),\qquad
\ggamma = (\alpha_1,\,\alpha_2,\,\alpha_3-1),\qquad\text{and}\qquad
\ddelta = (\alpha_1,\,\alpha_2-1,\,\alpha_3)\,.
$$
\label{prop:xDec}
\end{prop}

\proof Let us compute the gradient of $m_{\bbeta}$
\begin{equation*}
\nabla m_{\bbeta} = \frac{1}{h_P}\left(
\begin{array}{r}
(\alpha_1+1)\,m_{\aalpha} \\
\alpha_2\,m_{\ddelta_1} \\ 
\alpha_3\,m_{\ggamma_1}
\end{array}\right)\,,
\end{equation*}
and the cross products
\begin{equation*}
\quad
\mId \wedge
\left(
\begin{array}{c}
0 \\
m_{\ggamma} \\ 
0
\end{array}\right) =
\left(
\begin{array}{c}
-m_{\aalpha_{\phantom{1}}} \\
0 \\ 
\phantom{-}m_{\ggamma_1}
\end{array}\right)		
\quad\text{and}\quad
\mId \wedge
\left(
\begin{array}{c}
0 \\
0 \\ 
m_{\ddelta}
\end{array}\right) =
\left(
\begin{array}{c}
\phantom{-}m_{\aalpha_{\phantom{1}}} \\
-m_{\ddelta_1} \\ 
 0
\end{array}\right)\,,
\end{equation*}
where we defined the multi-indexes
$$
\ggamma_1 = (\alpha_1+1,\,\alpha_2,\,\alpha_3-1)\qquad\text{and}\qquad
\ddelta_1 = (\alpha_1+1,\,\alpha_2-1,\,\alpha_3)\,.
$$
Then to get Equation~\eqref{eqn:decForX}, 
we make a linear combination of such polynomial vectors
$$
c_1\,\frac{1}{h_P}\left(
\begin{array}{r}
(\alpha_1+1)\,m_{\aalpha} \\
\alpha_2\,m_{\ddelta_1} \\ 
\alpha_3\,m_{\ggamma_1}
\end{array}\right)
+c_2\,
\left(
\begin{array}{c}
-m_{\aalpha_{\phantom{1}}} \\
0 \\ 
\phantom{-}m_{\ggamma_1}
\end{array}\right)
+c_3\,
\left(
\begin{array}{c}
\phantom{-}m_{\aalpha_{\phantom{1}}} \\
-m_{\ddelta_1} \\ 
 0
\end{array}\right)\,.\\
$$
We observe that the particular choice of the multi-indexes 
$\bbeta,\,\ggamma$ and 
$\ddelta$,
leads to a linear combination of vectorial monomials 
which have monomials with the same multi-index on each component.
To complete the proof,
we solve the following linear system in the variables $c_1,\,c_2$ and $c_3$
\begin{equation*}
\left\{
\begin{array}{rl}
\mathlarger{\frac{1}{h_P}}(\alpha_1+1)\,c_1 - c_2 + c_3 =& 1 \\[1em]
\mathlarger{\frac{1}{h_P}}\alpha_2\,c_1 - c_3           =& 0 \\[1em]
\mathlarger{\frac{1}{h_P}}\alpha_3\,c_1 + c_2           =& 0 
\end{array}
\right.\,.
\end{equation*}
\cvd

In Propositions~\ref{prop:yDec} and~\ref{prop:zDec}, we provide 
similar results for vectorial monomials which have the other components different from 0.
We do not show the proofs of such propositions 
since they are similar to the one of Proposition~\ref{prop:xDec}.

%----------------------------------------------
%            SECONDA COMPONENTE
%----------------------------------------------
\begin{prop}
Considering a vectorial scaled monomial with only the \emph{second} component different from 0,
Equation~\eqref{eqn:decomp} becomes
$$
\left(
\begin{array}{c}
0 \\
m_{\aalpha} \\ 
0 
\end{array}\right) = 
\left(\frac{h_P}{|\aalpha|+1}\right) \nabla m_{\bbeta} + \left(\frac{\alpha_3}{|\aalpha|+1}\right) \mId \wedge
\left(
\begin{array}{c}
m_{\ggamma} \\
0 \\ 
0
\end{array}\right) - \left(\frac{\alpha_1}{|\aalpha|+1}\right) \mId \wedge
\left(
\begin{array}{c}
0 \\
0 \\ 
m_{\ddelta}
\end{array}\right)
$$
where 
$$
\bbeta  = (\alpha_1,\,\alpha_2+1,\,\alpha_3),\qquad
\ggamma = (\alpha_1,\,\alpha_2,\,\alpha_3-1),\qquad\text{and}\qquad
\ddelta = (\alpha_1-1,\,\alpha_2,\,\alpha_3)\,.
$$
\label{prop:yDec}
\end{prop}

%----------------------------------------------
%              TERZA COMPONENTE
%----------------------------------------------
\begin{prop}
Considering a vectorial scaled monomial with only the \emph{third} component different from 0,
Equation~\eqref{eqn:decomp} becomes
$$
\left(
\begin{array}{c}
0 \\
0 \\ 
m_{\aalpha} 
\end{array}\right) = 
\left(\frac{h_P}{|\aalpha|+1}\right) \nabla m_{\bbeta} + \left(\frac{\alpha_1}{|\aalpha|+1}\right) \mId \wedge
\left(
\begin{array}{c}
0 \\
m_{\ggamma} \\ 
0
\end{array}\right) - \left(\frac{\alpha_2}{|\aalpha|+1}\right) \mId \wedge
\left(
\begin{array}{c}
m_{\ddelta} \\
0 \\ 
0
\end{array}\right)
$$
where 
$$
\bbeta = (\alpha_1,\,\alpha_2,\,\alpha_3+1),\qquad
\ggamma = (\alpha_1-1,\,\alpha_2,\,\alpha_3),\qquad\text{and}\qquad
\ddelta = (\alpha_1,\,\alpha_2-1,\,\alpha_3)\,.
$$
\label{prop:zDec}
\end{prop}

\begin{rem}
A decomposition of a vectorial scaled monomial with more than one component different from zero can be obtained by summing 
the decompositions provided by the Propositions~\ref{prop:xDec},~\ref{prop:yDec} and~\ref{prop:zDec}.
\end{rem}

\begin{rem}
Since the first term of decompositions~\ref{prop:xDec},~\ref{prop:yDec} and~\ref{prop:zDec} are gradients,
we are able to generate \emph{any} vectorial polynomial $\p_{k}\in\GPerp$
starting from a linear combinations of $\mId\wedge\m$ where $\m\in[\mathcal{M}_{k-1}(P)]^3$.
\label{rem:GPerpPKm1}
\end{rem}

\subsection{Basis for $\mathcal{G}_k(P)$ and $\GPerp$}\label{sub:polyBasis}

The space relation provided in Equation~\eqref{eqn:spaceDec} suggests another polynomial vectorial basis for $[\Poly_k(P)]^3$.
Indeed, one can think to combine the basis of $\G_k(P)$ and $\GPerp$ 
to get a new basis of $[\Poly_k(P)]^3$, i.e.
$$
\p_k = \sum_{i=1}^{n} c_i\,\a_i + \sum_{j=1}^{m} d_i\,\b_i\,, 
$$
where $\{\a_i\}_{i=1}^{n}$ and $\{\b_j\}_{j=1}^{m}$ are two sets of vectorial basis function of 
$\G_k(P)$ and $\GPerp$, respectively.

\paragraph{Basis for $\mathcal{G}_k(P)$}
Since the operator $\nabla$ is an isomorphism between $\Poly_{k+1}(P)\backslash\R$ and $\G_k(P)$,
we know that a basis of $\Poly_{k+1}(P)\backslash\R$ will be mapped into a basis of $\G_k(P)$ by the operator $\nabla$.
Consequently, we have 
\begin{equation}
\G_k(P)=\text{span}\left\{\nabla m_{\aalpha}\right\}\,. 
\label{eqn:GBasis}
\end{equation}
where $0<|\aalpha|\leq k+1$.

\paragraph{Basis for $\GPerp$}
Starting from Equation~\eqref{eqn:decomp} and Remark~\ref{rem:GPerpPKm1},
the idea will be to exploit the basis of $[\Poly_{k-1}(P)]^3$ 
to get s basis for $\GPerp$.
However, there are some difficulties which require additional observations and resutls.

\begin{prop}
The linear function $\x\wedge *:[\Poly_{k-1}(P)]^3\to\GPerp$ is not an isomorphism.
\label{prop:noIso}
\end{prop}

\proof Suppose that $\x\wedge *$ is an isomorphism,
then $[\Poly_{k-1}(P)]^3$ and $\GPerp$ are isomorph and 
the following relation holds
\begin{equation*}
\dim([\Poly_{k-1}(P)]^3) = \dim(\GPerp)\,. 
\end{equation*}
Let us compute the dimensions of these two spaces separately
\begin{equation*}
\dim([\Poly_{k-1}(P)]^3) = \frac{k(k+1)(k+2)}{2} = \frac{k^3+4k^2+2k}{2}\,.
\end{equation*}
From Equations~\eqref{eqn:spaceDec} and~\eqref{eqn:GBasis} we have 
\begin{eqnarray}
\dim(\GPerp) &=& \dim([\Poly_k(P)]^3) - \dim(\G_k(P)) \nonumber \\[1em] &=& 
\frac{(k+1)(k+2)(k+3)}{2} - \left(\frac{(k+1)(k+2)(k+3)}{6} - 1 \right) \nonumber \\ &=& 
\frac{2k^3+9k^2+7k}{6}\,.
\label{eqn:dimGperp}
\end{eqnarray}
Since we do \emph{not} get the same dimension of the spaces,
$$
\frac{2k^3+9k^2+7k}{6} \neq \frac{k^3+4k^2+2k}{2}\,,
$$
$[\Poly_{k-1}(P)]^3$ and $\GPerp$ can not be isomorph and 
consequently $x\wedge *$ can not be an isomorphism.
\cvd

\begin{coroll}
Since $\x\wedge *$ is not an isomorphism, its kernel is not trivial and 
it is the set
\begin{equation}
\mathcal{K}_{k-1}(P):=\left\{ \p_{k-1}\in [\Poly_{k-1}(P)]^3\::\:  \p_{k-1} = \left(\begin{array}{c}
		 							    x\\
									    y\\
									    z
									  \end{array}\right)
									  \,p_{k-2}\,,
									  \forall p_{k-2}\in \Poly_{k-2}(P)\right\}\,.
\label{eqn:kerXCross}
\end{equation}
\label{corol:kerXCross}
\end{coroll}

\proof First of all we check if a generic element of $\mathcal{K}_{k-1}(P)$
is mapped to the null polynomial via the operator $x\wedge *$.
Let us consider a generic polynomial $p_{k-2}\in\Poly_{k-2}(P)$
$$
\left[\left(\begin{array}{c} x\\ y\\ z \end{array}\right)  \wedge \left(\begin{array}{c} x\\ y\\ z \end{array}\right) \,p_{k-2}\right] =
p_{k-2}\,\left[\left(\begin{array}{c} x\\ y\\ z \end{array}\right)  \wedge \left(\begin{array}{c} x\\ y\\ z \end{array}\right)\right] = 0
$$
Now we have to verify that only these elements are mapped to the null polynomial.
To achieve this goal,
we prove that $\mathcal{K}_{k-1}(P)$ has the same dimension of $\ker(x\wedge *)$.
Indeed, the dimension of $\mathcal{K}_{k-1}(P)$ is 
$$
\frac{k(k-1)(k+1)}{6}\,,
$$
and that the dimension of the $\ker(x\wedge *)$ is given by the following relation
$$
\dim{([\Poly_{k-1}(P)]^3)} = \dim(\GPerp) + \dim(\ker(x\wedge *))\,,
$$
which implies that 
\begin{eqnarray*}
\dim(\ker(x\wedge *)) &=& \dim{[\Poly_{k-1}(P)]^3} - \dim(\GPerp) \\[1em] &=&
\frac{k^3+3k^2+2k}{2} - \frac{2k^3+9k^2+7k}{6} \\ &=&
\frac{k(k-1)(k+1)}{6}\,.
\end{eqnarray*}
\cvd

\begin{rem}
Proposition~\ref{prop:noIso} still holds if we consider 
the linear operator $\mId\wedge *$ instead of $\x\wedge *$.
\end{rem}

Since $\mId\wedge *:[\Poly_{k-1}(P)]^3\to\GPerp$ is not an isomorphism,
a basis in the space $[\Poly_{k-1}(P)]^3$ is not mapped to a basis of $\GPerp$ so
we can not proceed in a similar way as for finding a basis for $\mathcal{G}_k(P)$.

\begin{prop}
Let us consider the a multi-index $\aalpha$, then the following relation holds
\begin{equation}
\mId\wedge \left(\begin{array}{c} m_{\aalpha} \\ 0          \\ 0           \end{array}\right) = -
\mId\wedge \left(\begin{array}{c} 0           \\ m_{\bbeta} \\ 0           \end{array}\right) - 
\mId\wedge \left(\begin{array}{c} 0           \\ 0          \\ m_{\ggamma} \end{array}\right) \,, 
\end{equation}
when 
$$
\bbeta = (\alpha_1-1,\alpha_2+1,\alpha_3)\qquad\text{and}\qquad\ggamma = (\alpha_1-1,\alpha_2,\alpha_3+1)\,.
$$
\label{prop:inTermsOfGPerpBasis}
\end{prop}

\proof 
The proof of this result is a simple computation of cross products.
\cvd

Proposition~\ref{prop:inTermsOfGPerpBasis} suggests us that we can replace the vectorial polynomials
$(m_{\aalpha},\,0,\,0)^T$ with a linear combination of other two vectorial monomials 
$(0,\,m_{\bbeta},\,0)^T$ and $(0,\,0,\,m_{\ggamma})^T$ in the decompositions 
of Propositions~\ref{prop:yDec} and~\ref{prop:zDec}.
Moreover, if we define the set of vectorial monomials
$$
\mathcal{M}^{\mathcal{G},\oplus}_k(P) := \left\{
\left(\hspace{-0.5em}\begin{array}{c} m_{\aalpha} \\ 0             \\ 0             \end{array}\hspace{-0.5em}\right)\,:\,
m_{\aalpha}\in\mathcal{N}_k(P)
\right\}\cup
\left\{
\left(\hspace{-0.5em}\begin{array}{c} 0           \\ m_{\bbeta}    \\ 0             \end{array}
\hspace{-0.5em}\right)\,:\,
m_{\bbeta}\in\mathcal{M}_k(P)
\right\}\cup
\left\{
\left(\hspace{-0.5em}\begin{array}{c} 0           \\ 0             \\ m_{\ggamma}   \end{array}
\hspace{-0.5em}\right)\,:\,
m_{\ggamma}\in\mathcal{M}_k(P)
\right\}\,,
$$
where
\begin{equation}
\mathcal{N}_k(P) :=\left\{m_{\aalpha}\::\: 0<|\aalpha|\leq k, \alpha_1 = 0\right\}\,,
\label{eqn:NDef}
\end{equation}
it follows that \emph{any} vector in $\GPerp$ can be written as a linear combination 
of the elements in $\mathcal{M}^{\mathcal{G},\oplus}_{k-1}(P)$, i.e. 
$\forall \p_k\in\GPerp$ we have 
$$
\p_k = c_0\,\mId\wedge\m_0^{\oplus} + c_1\,\mId\wedge\m_1^{\oplus} + \ldots + c_n\,\mId\wedge\m_n^{\oplus}\,,
$$
where $c_0,\,c_1,\ldots c_n\in\mathbb{R}$ and 
$\m_0^{\oplus},\,\m_1^{\oplus},\,\ldots\m_n^{\oplus}\in\mathcal{M}^{\mathcal{G},\oplus}_{k-1}(P)$.

Now if we show that the dimension of $\mathcal{M}^{\mathcal{G},\oplus}_{k-1}(P)$ coincides with 
the dimension of $\GPerp$, 
the image of the set $\mathcal{M}^{G,\oplus}_{k-1}(P)$ via the operator $\mId\wedge *$ is a basis for $\GPerp$.

\begin{prop}
The following space relation holds
\begin{equation}
\dim(\mathcal{M}^{\mathcal{G},\oplus}_{k-1}(P)) = \dim(\GPerp)\,.
\label{eqn:dimMAndG}
\end{equation}
\end{prop}
\proof We already show in the proof of Proposition~\ref{prop:noIso} that 
$$
\dim(\GPerp) = \frac{2k^3+9k^2+7k}{6}\,,
$$
then we have that 
$$
\dim(\mathcal{M}^{\mathcal{G},\oplus}_{k-1}(P)) = \frac{k(k+1)}{2} + 2\,\frac{k(k+1)(k+2)}{6} = 
\frac{2k^3+9k^2+7k}{6}\,,
$$
and this complete the proof.
\cvd

\section{Model problem: elliptic equation in mixed form}
\label{mod_prob}

The object of this work is the solution of the variational problem arising from the mixed formulation of a scalar elliptic equation in three spatial dimension. 

Let $\Omega$ be a bounded Lipschitz domain in $\RR^3$, whose boundary is denoted by $\partial\Omega$. We define the function spaces 
\[
\bV:=\{\u\in H(\div,\Omega)\::\:\,\,\u\cdot\n = u_N\quad\text{on }\partial\Omega\},
\]
and
\[
Q:=\left\{q\in L^2(\Omega)\::\: \int_\Omega q\,dx = 0\right\},
\]
where $H(\div,\Omega)$ is the space of vector-valued functions such that $\u$ and $\div{(\u)}$ belong to $[L^2(\Omega)]^3$ and $L^2(\Omega)$, respectively, and $u_N$ is the given Neumann datum.

The variational problem reads: 

\begin{equation}\label{mod_prob_eq}
\left\{
\begin{array}{rll}
\text{find } (\u,p)\in(\bV,Q):\\[0.1cm]
\displaystyle\a(\u,\v)-\b(\v,p) & \displaystyle= 0 & \forall\v\in\bV\vspace{0.2cm}\\
\displaystyle\b(\u,q) & \displaystyle=\int_\Omega f \, q \ dx & \forall q \in Q,
\end{array}
\right.
\end{equation}
where
\begin{equation}\label{bilin_form}
\begin{array}{l}
\displaystyle \a(\u,\v):=\int_\Omega \nu(x)\,\u\cdot\v\ dx,\vspace{0.2cm}\\
\displaystyle \b(\u,q):=\int_\Omega \div{(\u)} \, q\ dx,
\end{array}
\end{equation}
$f$ is a given function and $\nu(x)$ is a positive piecewise constant scalar function.

From the applications point of view, problem~\eqref{mod_prob_eq} arises in the context of multiphase incompressible flow through porous media, see e.g.~\cite{Bo-Bre-For}. Functions $\u$ and $p$ are usually called {\em velocity} and {\em pressure}. We refer to~\cite{Bo-Bre-For} for the mathematical analysis of such problem.

%-------------------------------------------------------------
\section{Virtual element discretization}\label{vem_disc}

Let $\Omega_h$ be a polyhedral decomposition of a three dimensional domain $\Omega$.
To solve Problem~\eqref{mod_prob_eq},
we follow a standard VEM approach.
We define local spaces in a generic polyhedron $P$, then 
we glue them together to get the global one.
Since we are dealing with the mixed formulation of the Laplace problem, 
we have to consider two types of spaces, one for the velocity $\VFG$ and one for the pressure $\VPG$.

\subsection{Definition of $\VFG$}\label{sub:velSpace}

To discretize the velocity, 
we take the virtual element 3d face space introduced in~\cite{divCurl}.
In this paper we will give a brief description on such local space, 
we refer to Section 5 of~\cite{divCurl} to have a deeper analysis.
However, since the definition of the degrees of freedom in~\cite{divCurl} is not appropriate from the 
practical and implementation point of view, 
we will make a more concrete definition of them.

Given a polyhedron $P$, we define the space
\begin{eqnarray}
\VFP:=\bigg\{\v_h\in H(\div;P)\cap H(\bcurl)&:& \v_h\cdot\n_f\in\Poly_k(f)\quad\forall f\in\partial P \nonumber \\
\phantom{\VFP:=\bigg\{\v_h\in H(\div;P)\cap H(\bcurl)}&\phantom{:}&
\div{(\v_h)}\in \Poly_{k-1}(P),\nonumber\\
\phantom{\VFP:=\bigg\{\v_h\in H(\div;P)\cap H(\bcurl)}&\phantom{:}& 
\bcurl(\v_h)\in [\Poly_{k-1}(P)]^3\bigg\}\,.\nonumber \\ 
\label{eqn:vectSpaceDef}
\end{eqnarray}

\noindent The degrees of freedom of such space are 
\begin{itemize}
\item normal face moments
 \begin{equation}
 \frac{1}{|f|}\mathlarger{\int_f} (\v_h\cdot\n_f)\,m_{\bbeta}\df\,,\qquad\forall f\in\partial P\,,\quad\forall m_{\bbeta}\in\mathcal{M}_k(f)\,,
 \label{eqn:facePoly}
 \end{equation}
 where $|f|$ denotes the area of the face $f$;
 \item internal gradient moments 
 \begin{equation}
 \frac{h_P}{|P|}\mathlarger{\int_P} \v_h\cdot\nabla m_{\aalpha}\dP\,,\qquad\forall m_{\aalpha}\in\mathcal{M}_{k-1}(P)\backslash\mathcal{M}_0(P)\,,
 \label{eqn:gradMom}
 \end{equation}
 where $|P|$ denotes the volume of $P$;
 \item internal cross moments 
 \begin{equation}
 \frac{1}{|P|}\mathlarger{\int_P} \v_h\cdot(\mathbf{m_I} \wedge \mathbf{m})\dP\,,\qquad\forall \mathbf{m}\in\mathcal{M}^{\mathcal{G},\oplus}_{k-1}(P)\,.
 \label{eqn:crossMom}
 \end{equation}
\end{itemize}

In Equations~\eqref{eqn:facePoly},~\eqref{eqn:gradMom} and~\eqref{eqn:crossMom},
we highlight the scaling factors.
In a virtual element framework the degrees of freedom have to scale as 1 
to get a better conditioning of the stiffness matrix~\cite{autostoppisti}.
Moreover, we explicitly show which polynomials are taken to define such degrees of freedom.
Since we have to consider a set of linearly independent conditions,
we use the basis functions provided in Subsection~\ref{sub:polyBasis}.
We will see that such choice makes computations easier 
and more straightforward with respect to the other ones.
This fact will become clearer when we show how to compute $\div(\v_h)$ and the projection operator $\bPi^0_k$.

\begin{rem}
The condition on the normal face moments could be replaced by the evaluation of $\v_h\cdot\n_f$ at suitable points on the face $f$~\cite{divCurl}.
However, finding a good position of such points could be not so straightforward when we are dealing with polygons so we use the degrees of freedom in Equation~\eqref{eqn:facePoly}.
\end{rem}

A generic function $\v_h\in\VFP$ is virtual so we can not use it directly.
To proceed with the VEM discretization of Problem~\eqref{mod_prob_eq}, 
we show that it is possible to compute some useful quantities.
\begin{itemize}
 \item \textbf{We can explicitly compute the polynomial} $(\v_h\cdot\n_f)$ 
 \textbf{for each face} $f$ \textbf{on} $\partial P$.\\  
 Since $(\v_h\cdot\n_f)\in\Poly_k(f)$, we can write it in terms of the monomial basis $\mathcal{M}_k(f)$.
 We exploit the normal face moments to find all coefficients $c_i$ of such polynomial, i.e. 
 we write such polynomial as
 $$
 (\v_h\cdot\n_f) = \sum_{|\ggamma|=0}^k c_i\,m_{\ggamma}^f
 $$
 and we test it against a each element of $\mathcal{M}_k(f)$
 $$
 \sum_{|\ggamma|=0}^k c_i \int_f m_{\ggamma}^f\,m_{\bbeta}^f\df = \int_f (\v_h\cdot\n_f)\,m_{\bbeta}^f\df\,,
 \quad\forall m_{\bbeta}^f\in\mathcal{M}_k(f)\,.
 $$
 Starting from these relations, we find \emph{exactly} $(\v_h\cdot\n_f)$ on each face of $\partial P$;
 \item \textbf{We can explicitly compute the divergence of a virtual function} $\v_h$.\\ 
 Since $\div(\v_h)~\in~\Poly_{k-1}(P)$, 
 we can write it in terms of the monomial basis $\mathcal{M}_{k-1}(P)$, i.e.
 $$
 \div(\v_h) = \sum_{|\ggamma|=0}^{k-1} c_i\,m_{\ggamma}\,,
 $$
 then, to get the coefficients $c_i$ of such polynomial, 
 we test it against each term of the scaled-monomial basis of degree $k-1$
 $$
 \sum_{|\ggamma|=0}^{k-1} c_i \int_P m_{\ggamma}\,m_{\aalpha}\dP = \int_P \div(\v_h)\,m_{\aalpha}\dP\,,
 \quad\forall m_{\aalpha}\in\mathcal{M}_{k-1}(P)\,.
 $$
 Even if $\v_h$ is virtual, it is possible to compute exactly the right hand sides starting from the 
 internal face moments and internal gradient moments.
 Indeed, if we integrate by parts we have
 $$
 \int_P \div(\v_h)\,m_{\aalpha}\dP = -\int_P \v_h\cdot\nabla m_{\aalpha} + \sum_{f\in\partial P}\int_f (\v_h\cdot\n_f)\,m_{\aalpha}\df\,,
 $$
 the first integral is an internal gradient moment and
 we can find the polynomial $(\v_h\cdot\n_f)$ from the internal face degrees of freedom,
 see the previous item.
 \item \textbf{We can compute an $L^2$-projection operator.}\\
 We define the $L^2$-projection operator $\bPi_k^0:\VFP\to[\Poly_k(P)]^3$ via 
 $$
 \int_P \bPi_k^0\,\v_h\cdot\p_k\dP = \int_P \v_h\cdot\p_k\dP\,,\qquad\forall \p_k\in\Poly_k(P)\,.
 $$
 To compute such projection operator, one considers the vectorial monomial base $[\M_k(P)]^3$ 
 for the projection $\bPi_k^0\,\v_h$, i.e.
 $$
 \bPi_k^0\,\v_h = \sum_{i=1}^{n_k} c_i\,\m_i\,,
 $$
 where $n_k=\dim([\M_k(P)]^3)$, and the relations 
 $$
 \sum_{i=1}^{n_k}\int_P \m_i \cdot \m_j\dP = \int_P \v_h\cdot\m_j\dP\,,\qquad\forall \m_j\in[\M_k(P)]^3\,,
 $$
 to find the unknown coefficients $c_i$.
 The right hand side of such conditions involves virtual function 
 so we have to understand if it is computable. 
 Exploiting Propositions~\ref{prop:xDec},~\ref{prop:yDec},~\ref{prop:zDec} and~\ref{prop:inTermsOfGPerpBasis}, we have 
 \begin{eqnarray}
 \int_P \v_h\cdot\m_j\dP &=& c_1\int_P \v_h\cdot\nabla m_{\bbeta}\dP 
 + c_2\int_P \v_h\cdot(\mathbf{m_I} \wedge \overline{\m})\dP\nonumber \\
 &+& c_3\int_P \v_h\cdot(\mathbf{m_I} \wedge \widetilde{\m})\dP\,,
 \label{eqn:decMonoForProj}
 \end{eqnarray} 
 where $c_1,c_2$ and $c_3$ are suitable constant,
 $m_{\bbeta}$, $\widetilde{\m}$ and $\overline{\m}$ are proper scaled monomials
 to decompose the vectorial monomial $\m_j$.
 The last two integrals in Equation~\eqref{eqn:decMonoForProj} are internal cross moments degrees of freedom.
 The first integral is a gradient moment only if $|\bbeta|\leq k-1$, otherwise we integrate by parts and get
 $$
 \int_P \v_h\cdot\nabla m_{\bbeta}\dP = -\int_P \div(\v_h)\,m_{\bbeta} + \sum_{f\in\partial P}\int_f(\v_h\cdot\n_f)\,m_{\bbeta}\df\,,
 $$
 which is still computable since we know both $\div(\v_h)$ and $(\v_h\cdot\n_f)$ on each face of $P$.
 \end{itemize}

\vspace{1em}
\noindent Then, the discrete velocity global space is defined by gluing such local spaces, i.e.
$$
\VFG := \left\{\v_h\in H^1(\div,\,\Omega)\::\:\v_h|_P\in\VFP\ \forall P\in\Omega_h,
\ \v_h\cdot\n = u_N,\text{ on }\partial\Omega\right\}\,.
$$

\subsection{Definition of $\VPG$}

To discretize the pressure, 
we consider a discontinuous polynomial space defined on each polyhedron $P$ of the discretization $\Omega_h$.
Given a polyhedron $P$, we introduce the local space
\begin{eqnarray}
\VPP := \left\{q_h\in L^2(P)\::\: q_h\in\Poly_{k-1}(P) \right\}\,.
\label{eqn:presSpace}
\end{eqnarray}
The degrees of freedom of such space are 
\begin{itemize}
 \item internal moments 
 \begin{equation}
 \frac{1}{|P|}\mathlarger{\int}_P q_h\,m_{\aalpha}\dP\qquad\forall m_{\aalpha}\in\mathcal{M}_{k-1}(P)\,.
 \label{eqn:presMom}
 \end{equation}
\end{itemize} 
Since the a function $q_h\in\VPP$ is a polynomial of degree $k-1$, 
it can be computed exactly on the polyhedron $P$ via the internal moments so 
we do not need to define any projection operator.

Then, as we have done for the global velocity discrete space,
the global pressure space is defined by gluing local spaces, i.e.
$$
\VPG := \left\{q_h\in L^2(\Omega)\::\: q_h|_P\in\VPP,\quad\forall P\in\Omega_h\text{ and }
\int_\Omega q_h dx = 0\right\}\,.
$$

\subsection{The discrete local forms}\label{sub:discForms}

To proceed with the discretization of Problem~\eqref{mod_prob_eq}, 
we construct suitable discrete forms~\cite{volley,autostoppisti,mixedVEM2d}.
As in standard VEM approach such discrete forms are defined element-wise and 
they depend on the degrees of freedom and projection operators.

Recalling Equation~\eqref{bilin_form}, we define 
\begin{eqnarray*}
\a_{h,P}(\v_h,\w_h) &:=& \nu\int_P \bPi_k^0\v_h\cdot\bPi_k^0\w_h\dP + 
s_P(\v_h -  \bPi_k^0\v_h, \w_h -  \bPi_k^0\w_h)\,,\\
\b_{h,P}(\v_h, q_h) &:=& \phantom{\nu}\int_P \div{(\v_h)}\,q_h\dP\,,\\
\f_{h,P}(q_h) &:=& \phantom{\nu}\int_P f\,q_h\dP\,,
\end{eqnarray*}
where $\v_h,\w_h\in\VFP,\,q_h\in\VPP$ and $s_P$ can be \emph{any} symmetric and positive definite bilinear form
which scales as the $\a_{P}(\cdot,\,\cdot)$. 
Such operator has to verify that 
there exist two constant $\alpha_*,\alpha^*>0$ such that
$$
\alpha_* a_P(\v_h,\,\v_h) \leq s_P(\v_h,\,\v_h) \leq \alpha^* a_P(\v_h,\,\v_h)\,,
$$
the coefficients $\alpha_*,\alpha^*$ depend on $\nu$ but not on the mesh-size.
In this paper we choose the Euclidean scalar product associated with the degrees of freedom of $\VFP$ multiplied by
the volume of $P$ and the value of $\nu$ at its barycenter~\cite{autostoppisti,mixedVEM2d}, i.e.
$$
s_P(\v_h,\,\w_h) := \nu(\x_P)\,|P|\,\sum_{i=1}^{\#dof_P} dof_i(\v_h)\,dof_i(\w_h)\,,
$$
where $\#dof_P$ are the number of degrees of freedom associated with a function in $\VFP$ and $dof_i:\VFP\to\mathbb{R}$ 
is a linear functional which associate to a function in $\VFP$ the value of its $i-$th degree of freedom.

\begin{rem}
The operators $\b_{h,P}$ and $\f_{h,P}$ involves polynomials ($\div{(\v_h)}$ and $q_h$) and the datum $f$.
Consequently, the discrete approximations of such terms 
are related only on the computation of integrals, i.e.
on the quadrature formulas used.
\end{rem}

\noindent Once we have defined the global forms
\begin{eqnarray*}
\a_h(\v_h,\w_h) &:=& \sum_{P\in\Omega_h} \a_{h,P}(\v_h,\w_h)\,,\\
\b_h(\v_h, q_h) &:=& \sum_{P\in\Omega_h} \b_{h,P}(\v_h, q_h) \,,\\
\f_h(q_h) &:=& \sum_{P\in\Omega_h} \f_{h,P}(q_h)\,,
\end{eqnarray*}
the discrete variational problem reads 
\begin{equation}\label{mod_prob_eqDisc}
\left\{
\begin{array}{rll}
\text{find } (\u_h,p_h)\in \VFG\times\VPG:&&\\[0.1cm]
\a_h(\u_h,\v_h)-\b_h(\v_h,p_h) & = 0 & \forall\v_h\in\VFG\vspace{0.2cm}\\
\b_h(\u_h,q_h) & =\,(\f_h\,,q_h) & \forall q_h \in \VPG\,.
\end{array}
\right.
\end{equation}

% \todo[inline]{MOLTO IMPORTANTE: sembra molto streminzito e si deve aggiungere un risultato sulla convergenza, 
% che devo chiedere a Lorenzo e Alessandro. MA in questo modo si vede che il focus è sulla definizione delle basi
% di monomi, lo sviluppo di quelle formule complicate per la decomposizione e il parallelismo}
% 
% \todo[inline]{Provarei ad aggiungere questo}
% 
% \begin{theo}
% Consider a polyhedral mesh $\Omega_h$ where all the elements are uniformly star shaped with respect to a ball and 
% all the edge/face diameters are comparable with respect to the polyhedron diameters,
% then Problem~\eqref{mod_prob_eqDisc} has a unique solution $(\v_h,\,q_h)\in \VFG\times \VPG$ and 
% the following error estimates hold
% \begin{eqnarray*}
% ||\u-\u_h||_0 &\leq& C\,h^k\,\left(||\u||_{k+1} + ||q||_{k+1} \right)\,,\\
% ||q-q_h||_0   &\leq& C\,h^{k+1}\,\left(||\u||_{k+1} + ||q||_{k+1} \right)\,,\\
% ||\div(\u-\u_h)||_0 &\leq& C\,h^{k+1}\,\left(|f|_{k+1} + ||q||_{k+1} \right)\,,\\
% \end{eqnarray*}
% where $C$ is a constant depending on $\nu$ but independent from the mesh size $h$.
% \label{rem:conv}
% \end{theo}
% 
% \todo[inline]{Se dovessimo dimotrare queste convergenze verrebbe un libro...
% possiamo riconduci a quelle fatte nello articolo
% ``Mixed Virtual Element Methods for general second order elliptic problem on polygonal meshes''?}
% 
% \todo[inline]{In alternativa abbiamo pensato a questo remark,
% ma ci manca una citazione che Lorenzo ha dato a Simone sul treno}

\begin{rem}
Consider a polyhedral mesh $\Omega_h$ where all the elements are uniformly star shaped with respect to a ball and 
all the edge/face diameters are comparable with respect to the polyhedron diameters.
Under these assumptions Problem~\eqref{mod_prob_eqDisc} 
has a unique solution $(\v_h,\,q_h)\in \VFG\times \VPG$ satisfying the following error estimates:
\begin{eqnarray*}
||\u-\u_h||_0 &\leq& C\,h^k\,\left(||\u||_{k+1} + ||q||_{k+1} \right)\,,\\
||q-q_h||_0   &\leq& C\,h^{k+1}\,\left(||\u||_{k+1} + ||q||_{k+1} \right)\,,\\
||\div(\u-\u_h)||_0 &\leq& C\,h^{k+1}\,\left(|f|_{k+1} + ||q||_{k+1} \right)\,,\\
\end{eqnarray*}
where $C$ is a constant depending on $\nu$ but independent from the mesh size $h$.
The proof is beyond the scope of this paper and 
it could be obtained by combining the results in~\cite{da2017virtual,brenner2018virtual}.
\label{rem:conv}
\end{rem}

\subsection{Numerical results}\label{conv_numerical_results}

\begin{figure}[!htb]
\def\sizeOfFigMeshes{0.35}
\begin{center}
\begin{tabular}{cc}
\includegraphics[width=\sizeOfFigMeshes\textwidth]{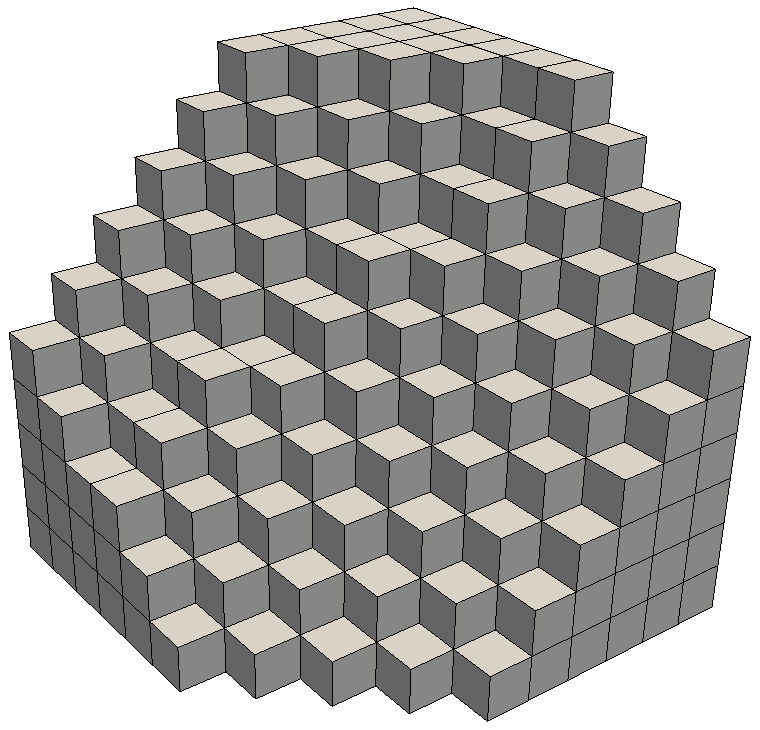}&
\includegraphics[width=\sizeOfFigMeshes\textwidth]{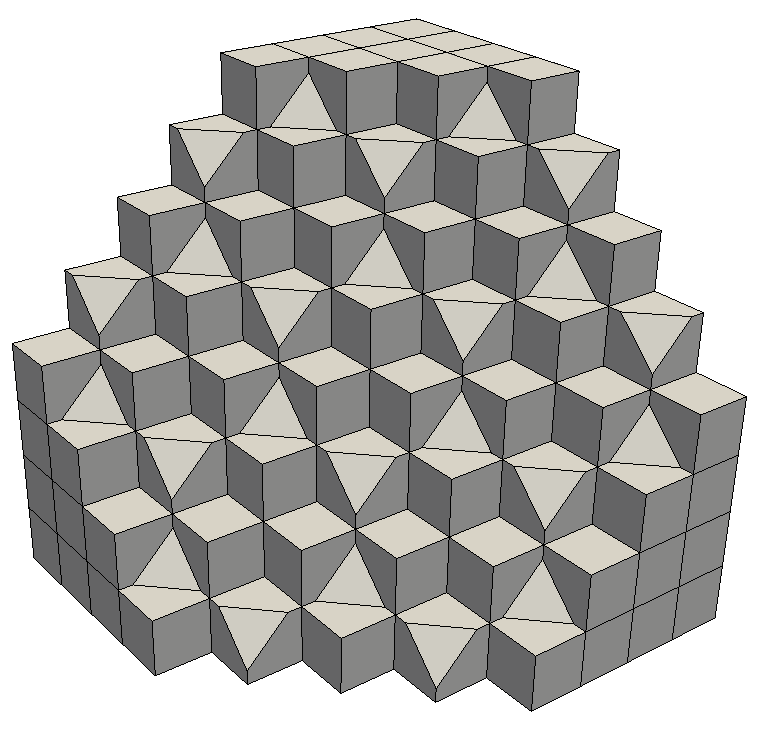}\\
\Cube{} &\Octa{}\\[1em]
\includegraphics[width=\sizeOfFigMeshes\textwidth]{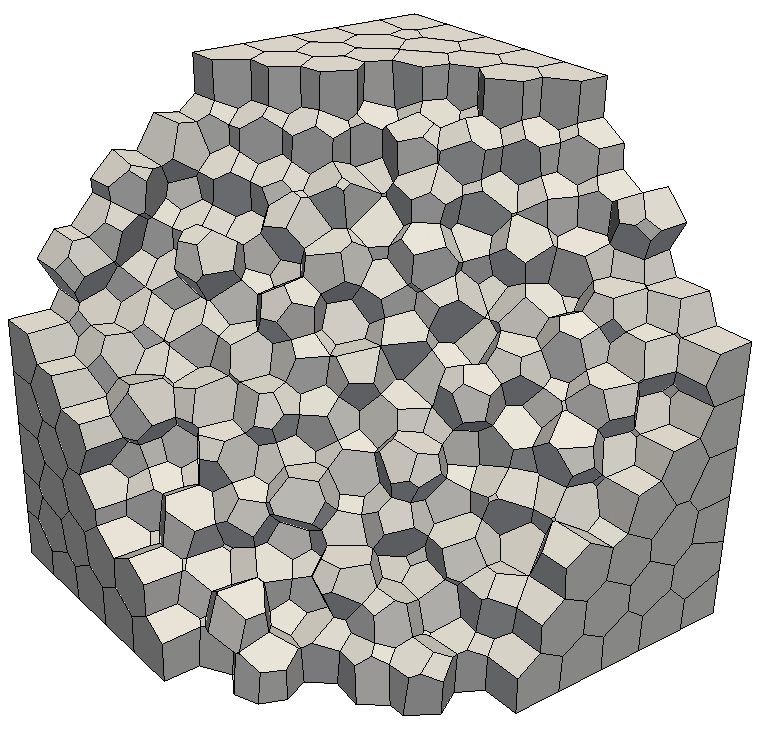}&
\includegraphics[width=\sizeOfFigMeshes\textwidth]{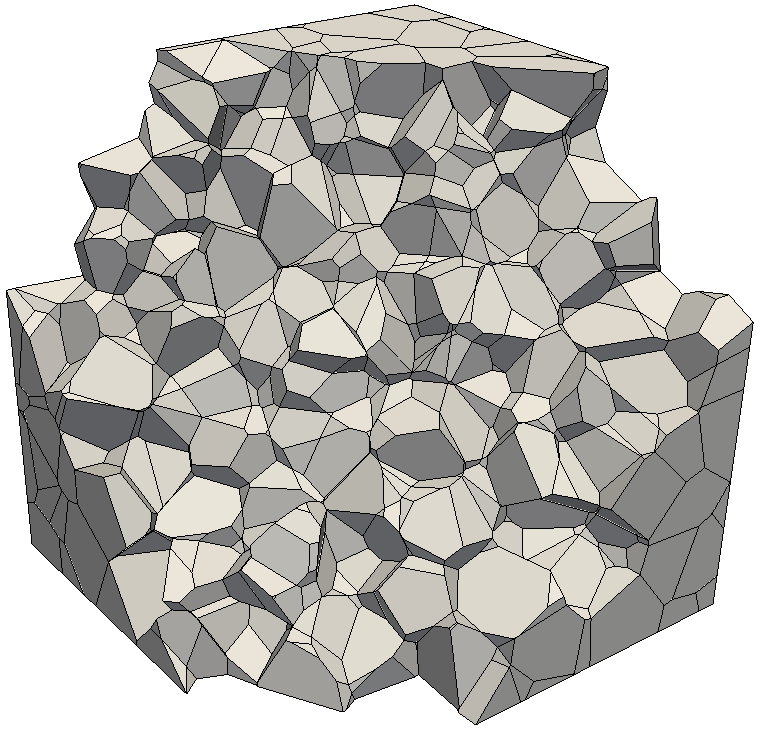}\\
\CVT{} &\Random{}
\end{tabular}
\caption{A sample of the mesh taken into account.}
\label{fig:meshes}
\end{center}
\end{figure}

We conclude this section with a numerical example 
to validate the mixed virtual element approach in solving Problem~\eqref{mod_prob_eq}.
In the following test we use four different discretizations of the cube $[0,\,1]^3$:
\begin{itemize}
 \item \Cube{}, a mesh composed by structured cubes,
 \item \Octa{}, a mesh composed by polyhedron with seven or eight faces,
 \item \CVT{}, a Voronoi tessellation optimized via a standard Llyod algorithm~\cite{cvtPaper},
 \item \Random{} a Voronoi tessellation of a set of points randomly put inside $\Omega$.
\end{itemize}
In such discretizations the mesh elements become more irregular. 
Indeed, firstly we take into account  standard cubes and regularly shaped polyhedrons, 
\Cube{} and \Octa{} meshes.
Then, we move to \CVT{} and \Random{} meshes which are characterized by small edges,
stretched and small faces, see Figure~\ref{fig:meshes}.

All the meshes taken into account were generated via the c++ library \texttt{voro++}~\cite{voroPlusPlus} 
and putting the control points of the Voronoi cells in a proper way~\cite{apollo11,cvtPaper}.
In order to study the error convergence rate, 
we generate a sequence of four progressive refinements composed by approximately 27, 125, 1000 and 8000 polyhedrons
and we associate with them a mesh-size
\begin{equation}
h := \frac{1}{N_P} \sum_{i=1}^{N_P} h_P\,,
\label{eqn:meshSize}
\end{equation}
where $N_P$ is the number of polyhedrons in the mesh.

We compute the error of both velocity and pressure. 
More specifically we consider the following $L^2$ error indicators
\begin{equation}
e_{\v}:=\frac{\sqrt{\sum_{i=1}^{N_P} ||\v-\bPi_k^0\v_h||_{0,P}^2}}{||\v||_{0,\Omega}}\,,\qquad
e_{q}:=\frac{\sqrt{\sum_{i=1}^{N_P} ||q-q_h||_{0,P}^2}}{||q||_{0,\Omega}}\,,\qquad
\end{equation}
where $||\cdot||_{0,\mathcal{D}}$ denotes the standard $L^2$ norm over a domain $\mathcal{D}$.

We remark that the function $\v_h$ is virtual so we exploit its projection to compute the error $e_{\v}$.
On the other hand the discrete function associated with the pressure is a polynomial of degree $k-1$ on each element and we do not use any projection operator to compute such error.
The expected convergence rates of $e_{\v}$ and $e_q$ are $O(h^k)$ and $O(h^{k+1})$,
respectively, see Remark~\ref{rem:conv}.

We consider the Problem~\eqref{mod_prob_eq} with $\nu(\x)=1$ and 
we set both right hand side and the boundary conditions in such a way that the exact solution
is the couple
$$
\v(x,\,y,\,z):=\left(\begin{array}{c}
               -5x^4-y^2z^3\\
               -24y^3-2xyz^3\\
               -27z^2-3xy^2z^2\\
               \end{array}\right)
$$
and 
$$               
q(x,\,y,\,z):=x^5+6y^4+9z^3+xy^2z^3\,.
$$

\begin{figure}[!htb]
\def\sizeOfFigAll{0.45}
\begin{center}
\begin{tabular}{cc}
\includegraphics[width=\sizeOfFigAll\textwidth]{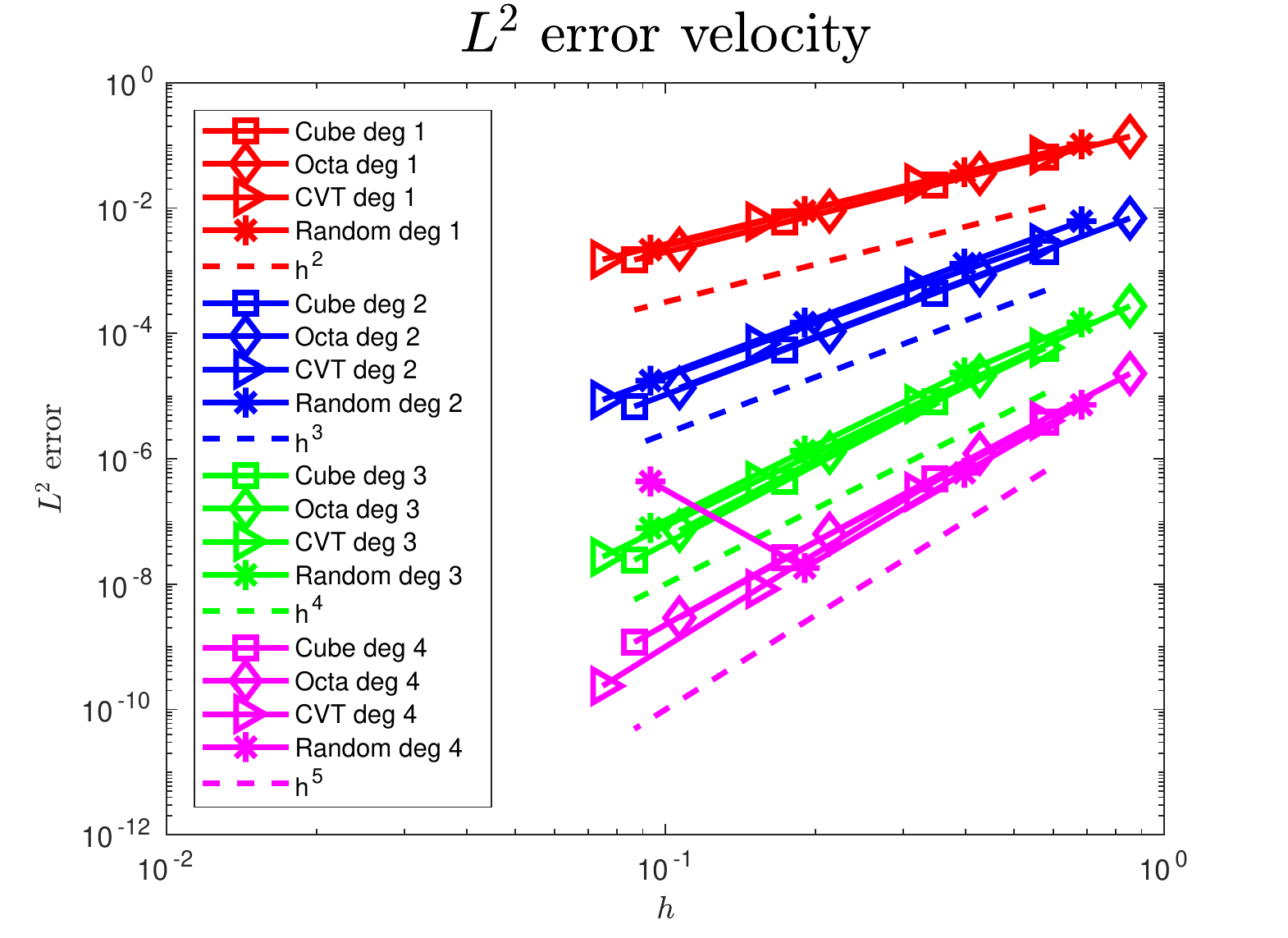}&
\includegraphics[width=\sizeOfFigAll\textwidth]{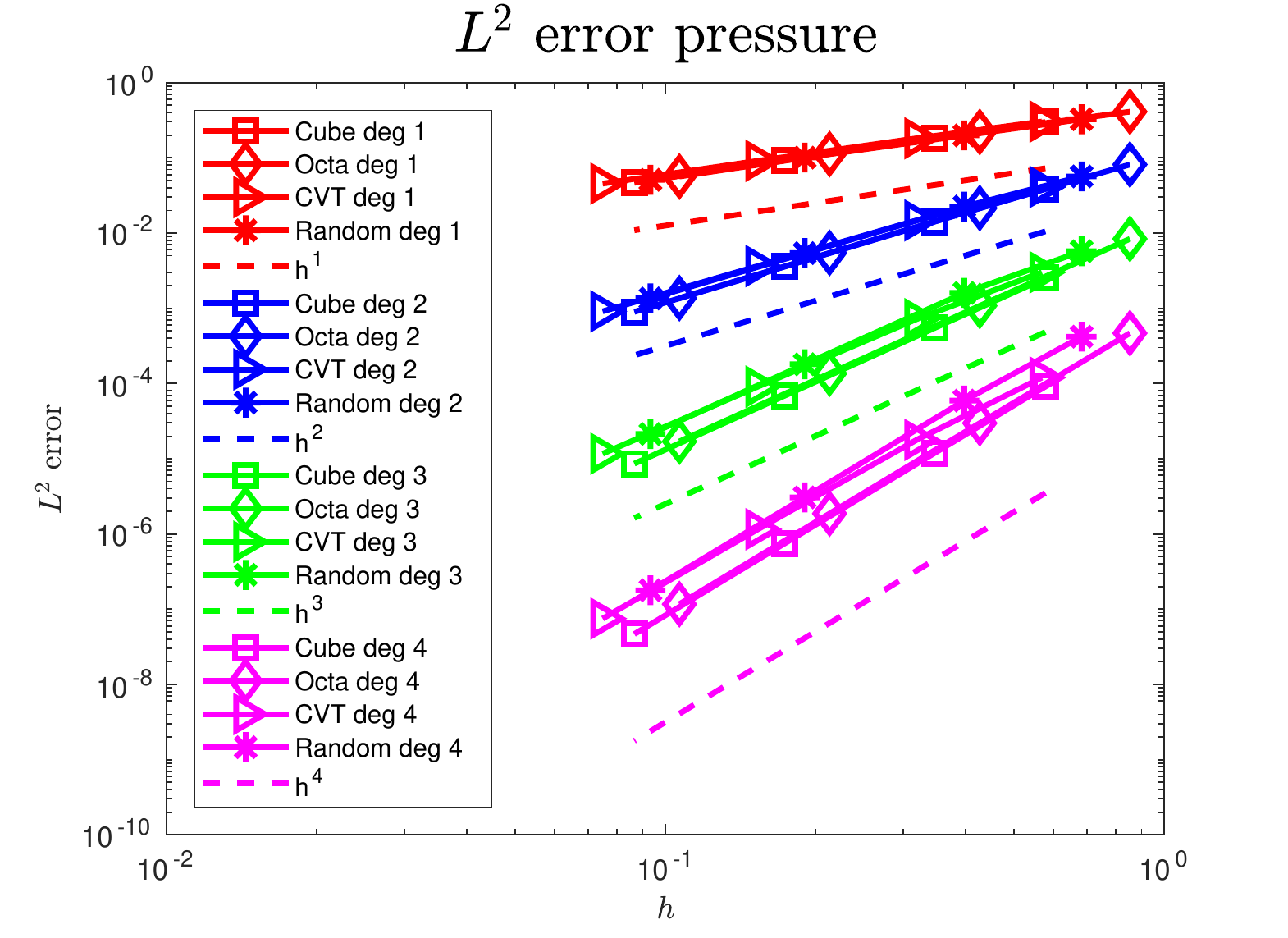}\\
\end{tabular}
\caption{Convergence lines for each set of meshes taken into account and degrees $k=2,3$ and 4.}
\label{fig:allConvAndDeg}
\end{center}
\end{figure}

In Figure~\ref{fig:allConvAndDeg} we show the convergence lines for all meshes and for $k=1,2,3$ and 4.
For all the set of meshes and for each approximation degree $k$,
the method behaves as expected and we recover the convergence rate predicted by the theory
for both velocity and pressure.

Moreover, if we fix the degree $k$ and we vary the type of meshes, 
these convergence lines are close to each other.
This fact is a further numerical prove of the robustness of VEM with respect to distorted elements.

The method fails only in the last step of the \Random{} meshes for the error $e_{\v}$.
This fact could be due to the ill-conditioning of the linear system at hand.
Indeed, such mesh is characterized by really small features 
(the smallest face has area $\num{3.9e-14}$ and the smallest edge is long $\num{5.9e-08}$)
which may affect the condition number of the stiffness matrix
and consequently the computation of the error,
when we have such an high approximation degree.

%-------------------------------------------------------------
\section{Parallel preconditioners}
\label{para_precon}

Our strategy for building an efficient parallel solver is based on the parallel library PETSc from Argonne National Laboratory~\cite{petsc-web-page,petsc-user-ref,petsc-efficient}.
Such library is built on the MPI standard and 
it offers advanced data structures and routines for the parallel solution of partial differential equations, from basic vector and matrix operations to more complex linear and nonlinear equation solvers. 
In our c++ code, vectors and matrices are built and subassembled in parallel on each processor.

Let us denote by $\mathcal{A}$ the linear system matrix arising from the discretization of the model problem (\ref{mod_prob_eq}), which has the typical saddle point structure
\[
\mathcal{A}=
\left[
\begin{array}{cc}
A & B^T\\
B & -C
\end{array}
\right].
\]
To solve such  linear system, we use the parallel GMRES method provided by the PETSc library,
preconditioned by two types of block-diagonal preconditioners (see e.g. \cite{benziGL.2005,mardal.2011,axelssonBBKA.2016}) of the form
\begin{equation}\label{B_D}
\mathcal{B}_D = \left[
\begin{array}{cc}
\mathrm{B}_{1} & 0\\
0 & \mathrm{B}_{2}\\
\end{array}
\right]:
\end{equation}
\begin{itemize}
\item {\bf Block-Schur} where
\begin{eqnarray}
B_1^{-1} & = & \mbox{ diagonal preconditioner for } A, \mbox{ i.e. }B_1=diag(A) \nonumber\\
B_2^{-1} & = & \mbox{ exact solution of the approximate Schur complement } S\label{block_schur}
\end{eqnarray}
with $S=-C-B~ diag(A)^{-1}B^T$. For the inversion of $S$ at each preconditioning step we use the parallel multifrontal direct solver Mumps~\cite{amestoy.2001,amestoy.2006}.
\item {\bf Block-Reg} where
\begin{eqnarray}
B_1^{-1} & = & \mbox{ Algebraic Multigrid preconditioner for } A+B^T W^{-1} B,\nonumber\\
B_2^{-1} & = & W^{-1}\label{block_diag}
\end{eqnarray}
with $W=\gamma I$, for a suitable parameter $\gamma>0$. As Algebraic Multigrid preconditioner we use the GAMG solver of PETSc.
\end{itemize}
In the following tests we compare the previous two block-diagonal preconditioners and the parallel direct solver Mumps 
considering the model problem in Subsection~\ref{conv_numerical_results}.

\subsection{Numerical results}
\label{par_numerical_results}

In the numerical tests we use the Linux cluster INDACO 
(www.indaco.unimi.it) of the University of Milan,
constituted by 16 nodes, each carrying 2 processors
INTEL XEON E5-2683 V4 2.1 GHz, with 16 cores each.

%\todo[inline]{Mi è sembrato opportuno anticipare questa parte visto che viene utilizzata in entrambi gli esempi.}
We consider three types of polyhedral meshes, {\bf Cube}, {\bf Octa} and {\bf CVT}, introduced in 
Subsection~\ref{conv_numerical_results}.
We solve our model problem~\eqref{mod_prob_eq} using the proposed VEM discretizations.
For the Block-Reg preconditioner, we heuristically found that 
the best performances are obtained taking $\gamma=h^2$, 
where $h$ is the mesh size parameter defined in Equation~\eqref{eqn:meshSize}.

\begin{table}[!htb]
\begin{center}
\begin{tabular}{c|c|cc|cc|ccc|ccc}
\hline
\multicolumn{12}{c}{{\bf Cube mesh with 32768 elements, k = 1, dofs = 435201}} \\
\hline
$p$ &$S_{p}^{\text{id}}$ &  $T_{ass}$  &$S_{p}$ & \multicolumn{2}{c|}{Mumps} & \multicolumn{3}{c|}{Block-Schur} &  \multicolumn{3}{c}{Block-Reg} \\%[1mm]
&  & & &$T_{sol}$  &$S_{p}$  & it & $T_{sol}$ & $S_{p}$ &it & $T_{sol}$ & $S_{p}$\\
\hline
1       &-      &68     &-      &823    &-      &66     &50     &-      &84     &53     &-\\
4       &4      &19     &3.6    &228    &3.6    &66     &17     &2.9    &114    &32     &1.7\\
8       &8      &10     &6.8    &138    &6.0    &66     &12     &4.2    &116    &29     &1.8\\
16      &16     &5      &13.6   &74     &11.1   &66     &10     &5.0    &120    &26     &2.0\\
32      &32     &3      &22.7   &47     &17.5   &66     &9      &5.5    &137    &26     &2.0\\
\hline
\hline
\multicolumn{12}{c}{{\bf Cube mesh with 13824 elements, k = 2, dofs = 508033}} \\
\hline
$p$ &$S_{p}^{\text{id}}$ &  $T_{ass}$  &$S_{p}$	& \multicolumn{2}{c|}{Mumps} & \multicolumn{3}{c|}{Block-Schur} &  \multicolumn{3}{c}{Block-Reg} \\%[1mm]
&  & & &$T_{sol}$  &$S_{p}$  & it & $T_{sol}$ & $S_{p}$ &it & $T_{sol}$ & $S_{p}$\\
\hline
1	&-       &271	&- 		&659	&-		&72	&347	&-	&98	&138	&-\\
4	&4	 &97	&2.8 	&436	&1.5 	&72	&118	&2.9 	&162	&73	&1.9 \\
8	&8	 &53	&5.1 	&270	&2.4 	&72	&77	&4.5 	&168	&48	&2.9 \\
16	&16	 &27	&10.0 	&157	&4.2 	&72	&55	&6.3  	&172	&67	&2.1 \\
32	&32	 &13	&20.8  	&95	&6.9 	&72	&40	&8.7 	&182	&58	&2.4 \\
%64	&64      &7	&38.7 (64) 	&205	&3.2 (64)	&72	&93	&3.7 (64) 	&796	&473	&1.4 (64)\\
\hline
\hline
\multicolumn{12}{c}{{\bf Cube mesh with 8000 elements, k = 3, dofs = 612001}} \\
\hline
$p$ &$S_{p}^{\text{id}}$ &  $T_{ass}$  & $S_{p}$ & \multicolumn{2}{c|}{Mumps} & \multicolumn{3}{c|}{Block-Schur} &  \multicolumn{3}{c}{Block-Reg} \\%[1mm]
&  & & &$T_{sol}$  &$S_{p}$  & it & $T_{sol}$ & $S_{p}$ &it & $T_{sol}$ & $S_{p}$\\
\hline
4	&-   &623	&-		&461	&-		&147	&570	&-	&720	&450	&-\\
8	&2   &334	&1.9 	&277	&1.7 	&147	&430	&1.3 	&784	&285	&1.6 \\
16	&4   &194	&3.2 	&163	&2.8 	&147	&268	&2.1 	&774	&203	&2.2 \\
32	&8   &97	&6.4 	&99	&4.7 	&147	&192	&3.0 	&918	&197	&2.3 \\
%64	&16  &49	&12.7 	&207	&2.2 	&147	&531	&1.1 	&980	&753	&0.6 \\
\hline
\end{tabular}
\caption{Strong scaling test, {\bf Cube} meshes. $p$:=number of procs; $S_{p}^{\text{id}}$:=ideal speedup; $T_{ass}$:=assembling time in seconds; $T_{sol}$:=solution time in seconds; it:=GMRES iterations; $S_p$:=parallel speedup computed with respect to the 1 procs run for k=1,2 and to the 4 procs run for k=3.}
\label{table_hexa_scal}
\end{center}
\end{table}

\begin{figure}[!htb]
\begin{center}
\includegraphics[width=0.45\textwidth]{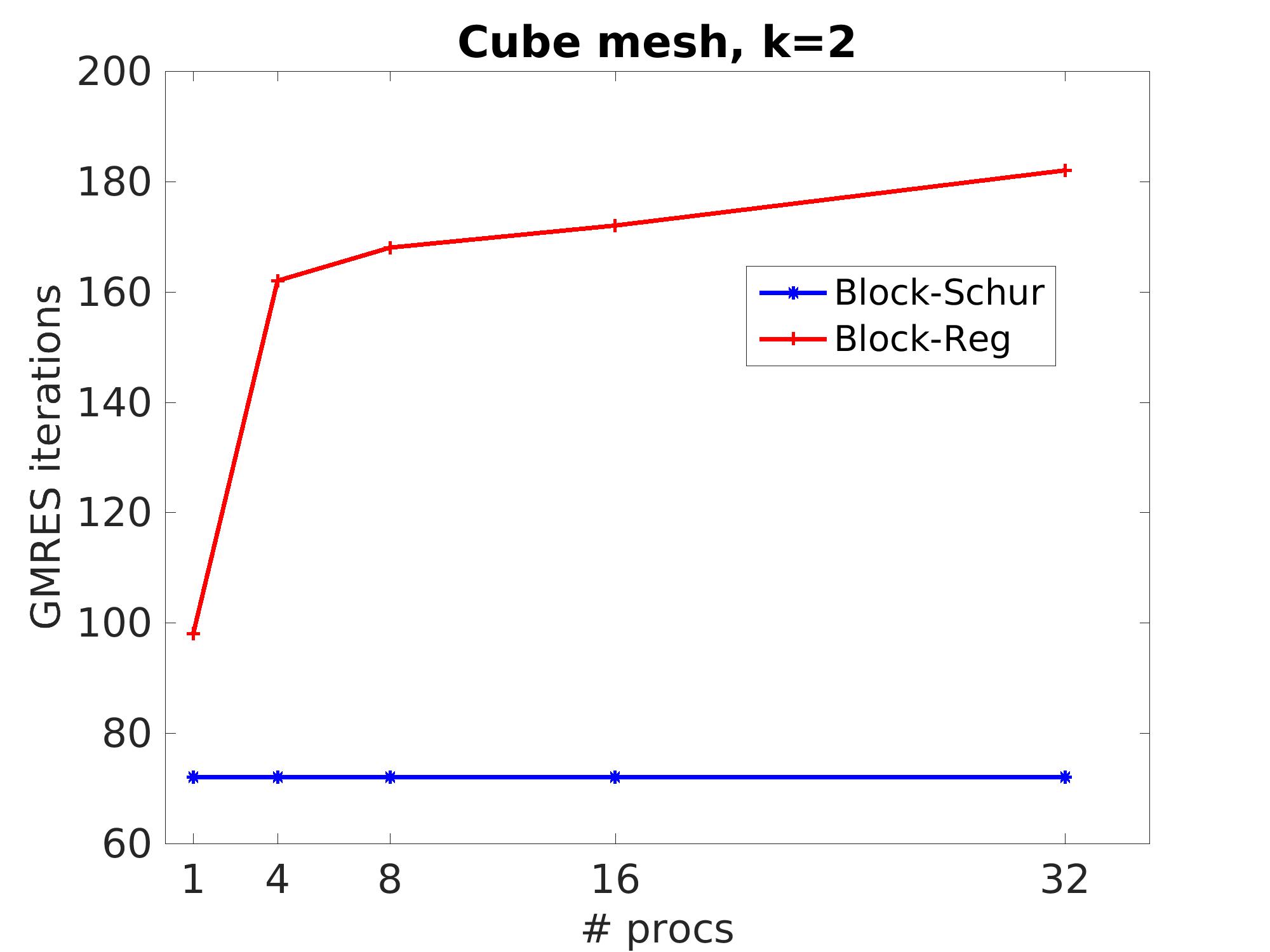}
\includegraphics[width=0.45\textwidth]{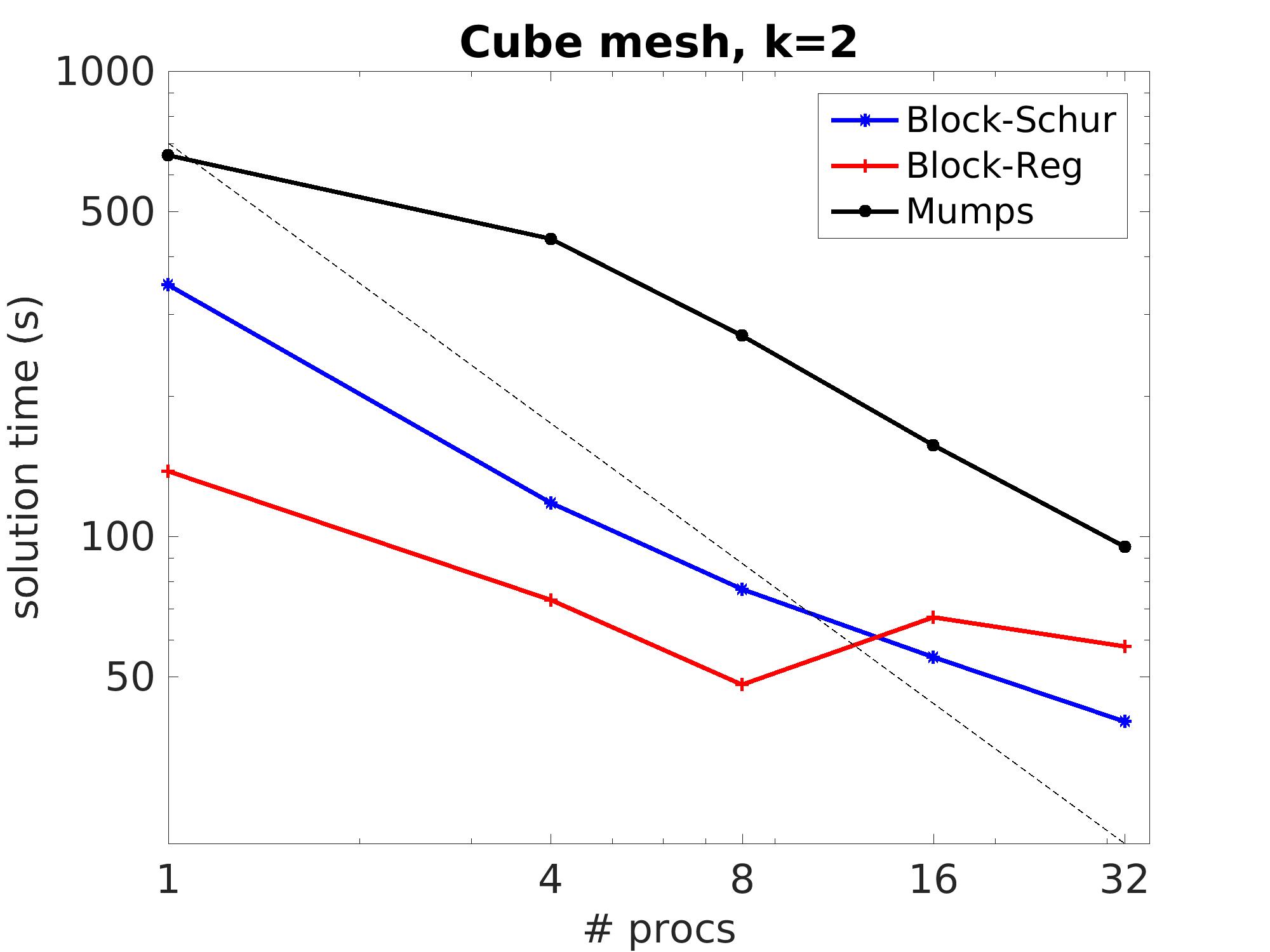}\\
\includegraphics[width=0.45\textwidth]{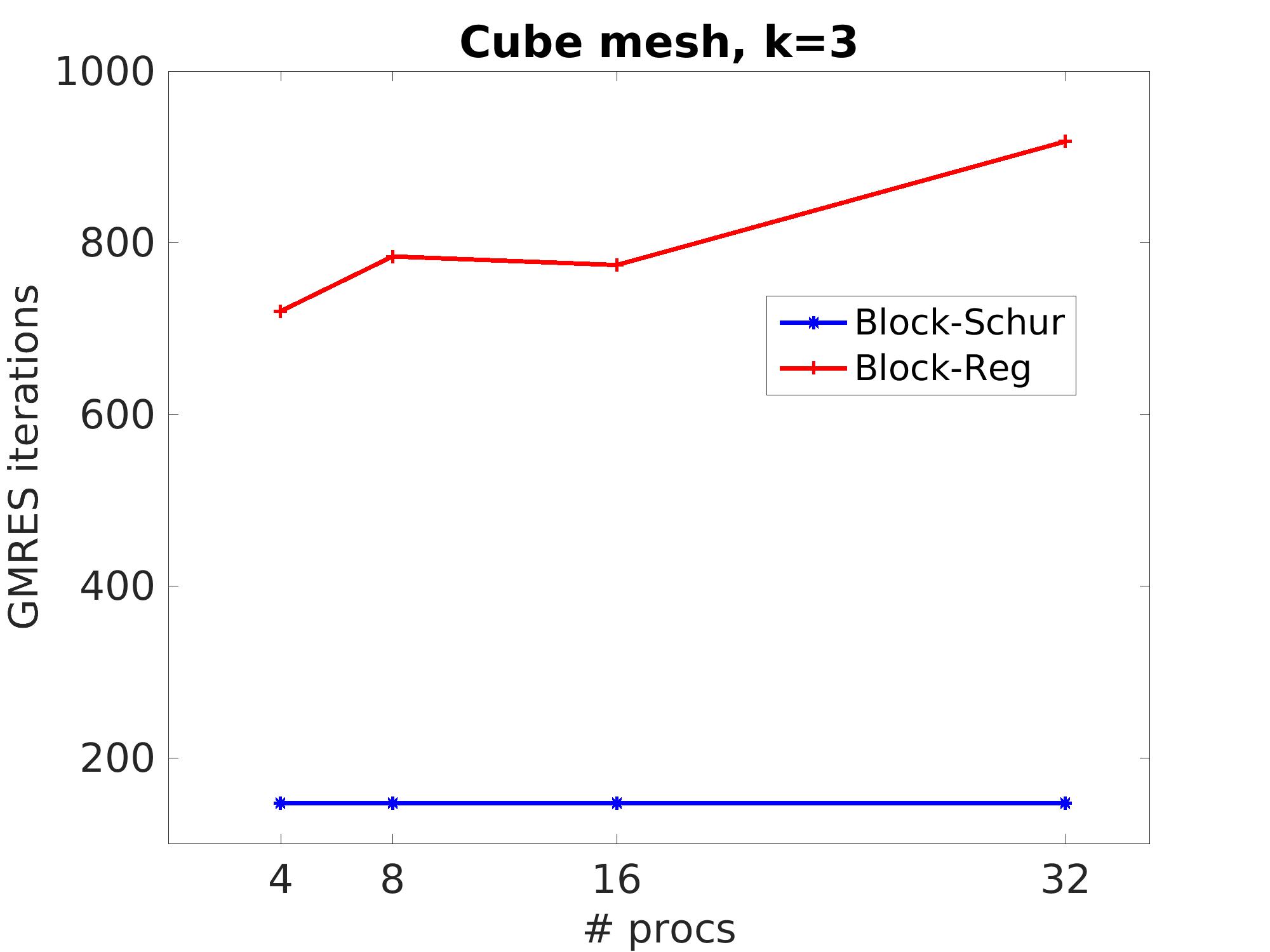}
\includegraphics[width=0.45\textwidth]{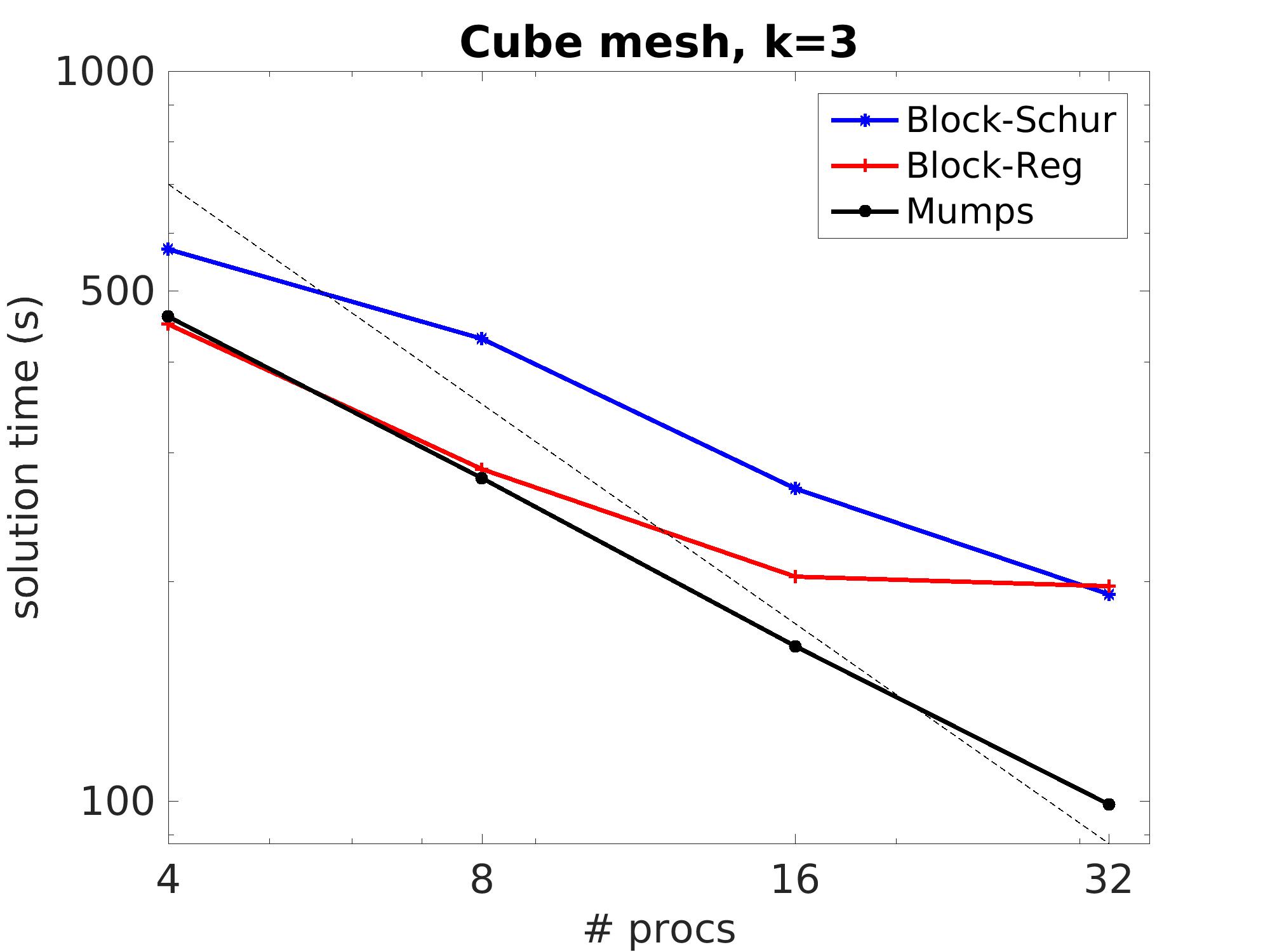}\\
\caption{Strong scaling test, {\bf Cube} meshes. GMRES iterations (left) and solution time (right) as a function of the number of procs for $k=2$ (first row) and $k=3$ (second row). In the solution time plots, the black dashed line indicates the steepness of the ideal time reduction.}
\label{fig_hexa_scal}
\end{center}
\end{figure}

\subsubsection{Test 1: strong scaling}

In this paragraph, we study the parallel performance of the three solvers (Mumps, Block-Schur and Block-Reg), by increasing the number of processors, while keeping fixed the global number of degrees of freedom (dofs). Hence, this is a {\em strong scaling} test. For all the three meshes, we consider a VEM discretizations of order $k=1,2,3$.

The results of the strong scaling test are displayed in Tables \ref{table_hexa_scal}, \ref{table_voro_scal} and \ref{table_octa_scal}. Note that, in case of order $k=3$, the runs with 1 processor (procs) went Out of Memory. Thus, we report the results starting from the 4 procs run. Denoting by $p$ the number of procs, we recall that the {\em parallel speedup} $S_p$ is defined as
\[
S_p := \frac{\mbox{CPU time with 1 procs}}{\mbox{CPU time with }p\mbox{ procs}}.
\]

We first observe that, irrespective of the kind of mesh and of the order $k$ of the VEM discretization, the CPU times needed to assemble the stiffness matrix and the right hand side ($T_{ass}$) are scalable, with good speedup values quite close to the ideal ones.

In case of {\bf Cube} meshes (Table \ref{table_hexa_scal} and Fig. \ref{fig_hexa_scal}), for $k=1,2$, the Block-Schur preconditioner is completely scalable in terms of GMRES iterations, which do not depend on $p$, whereas the Block-Reg preconditioner shows a slight increase of GMRES iterations. The solution times ($T_{sol}$) decrease with $p$ for all the three solvers, but the speedup values are far from the ideal ones. The highest speedup values are achieved by Mumps for $k=1$ and Block-Schur for $k=2$. In terms of CPU times, for both $k=1,2$, the most effective solver is Block-Schur, being between 2 and 16 times as fast as Mumps. For $k=3$ instead, the iterative solvers suffer in terms of GMRES iterations, due to the severe ill-conditioning of the stiffness matrix. Mumps is the most effective solver, 
being twice as fast as the iterative methods on 32 procs. 

\begin{figure}[!htb]
\begin{center}
\includegraphics[width=0.45\textwidth]{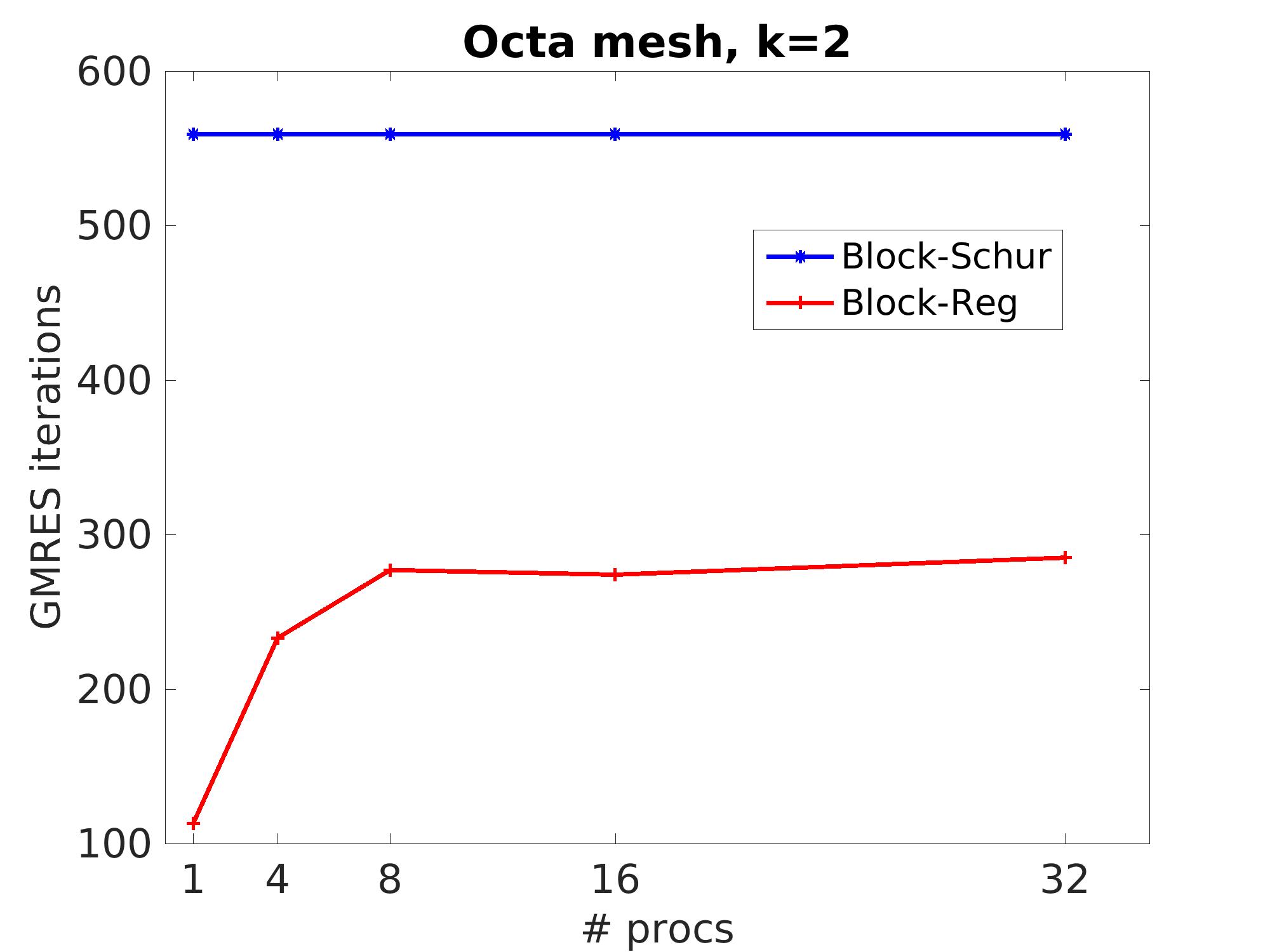}
\includegraphics[width=0.45\textwidth]{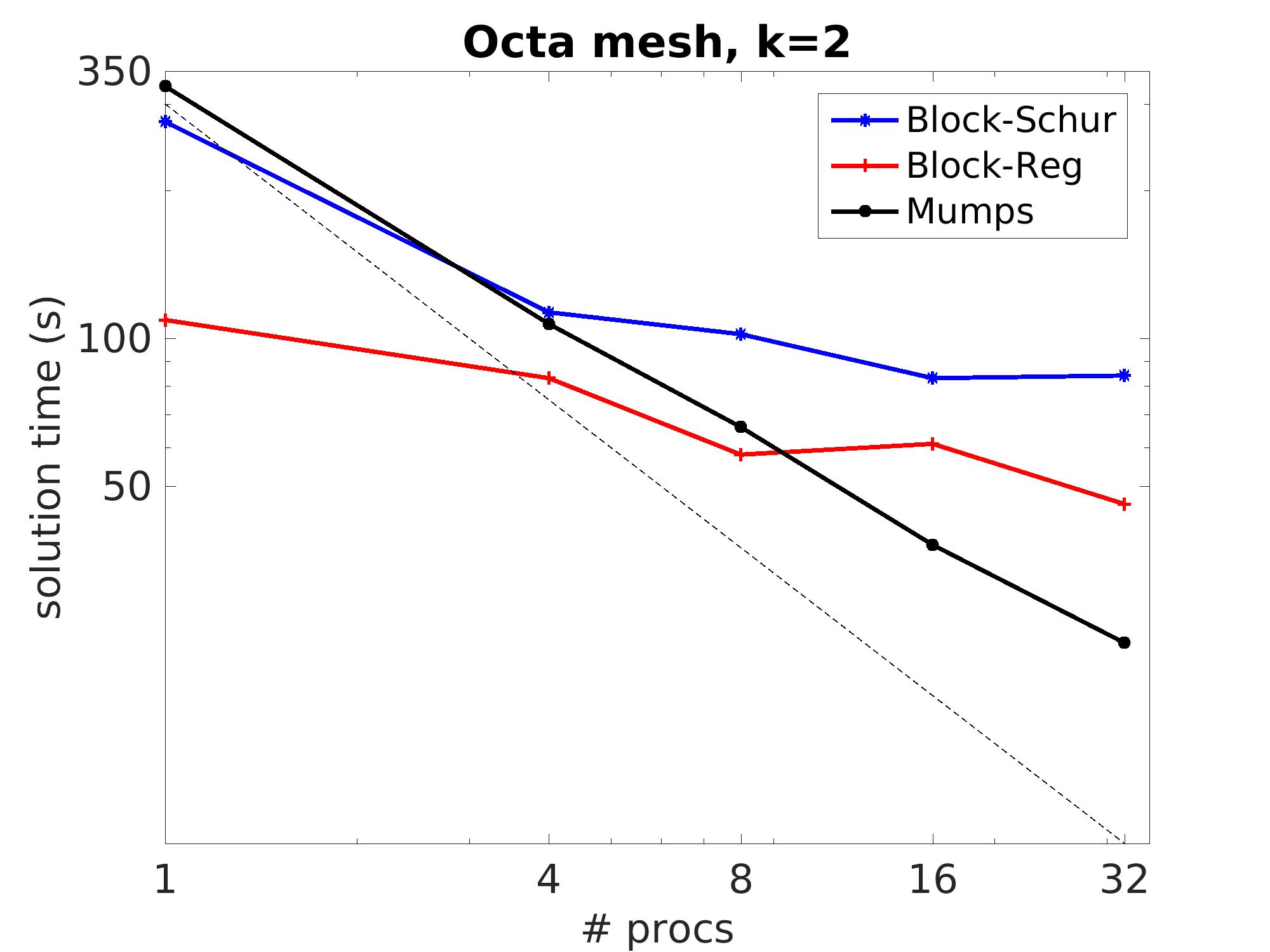}\\
\includegraphics[width=0.45\textwidth]{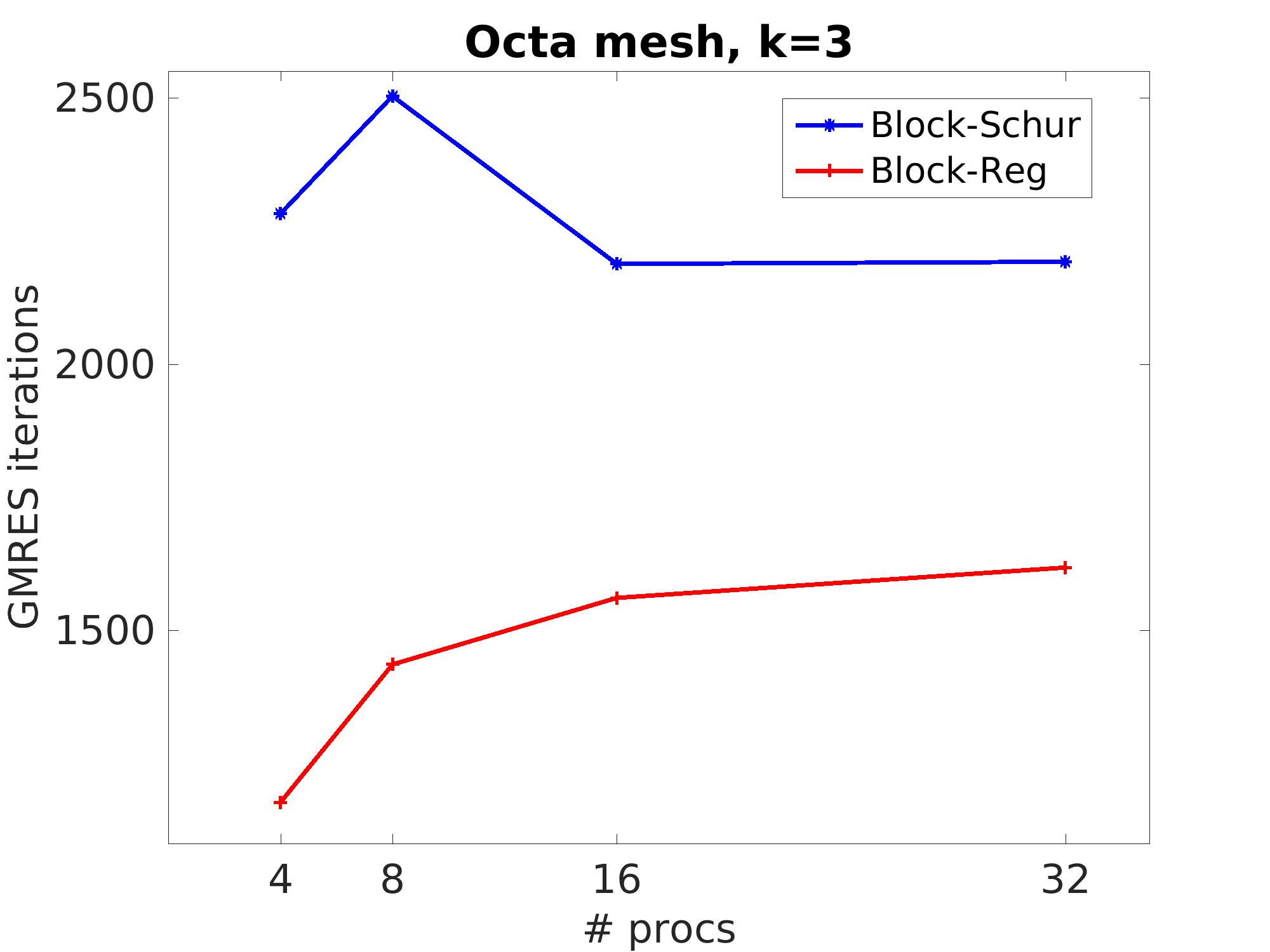}
\includegraphics[width=0.45\textwidth]{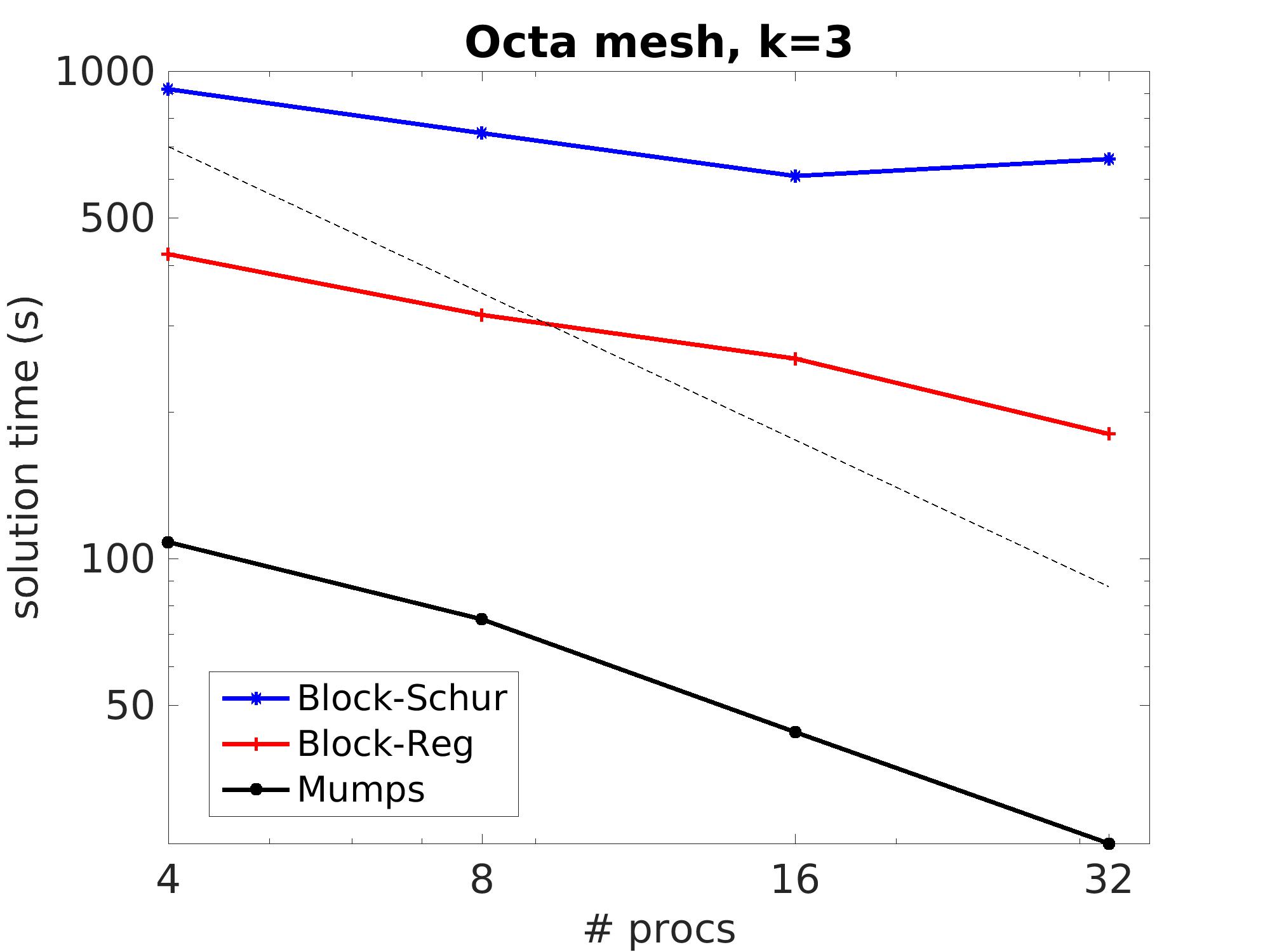}\\
\caption{Strong scaling test, {\bf Octa} meshes. Same format as in Fig. \ref{fig_hexa_scal}.}
\label{fig_octa_scal}
\end{center}
\end{figure}

In case of {\bf Octa} meshes (Table \ref{table_octa_scal} and Fig. \ref{fig_octa_scal}), 
the GMRES iterations of the iterative solvers are higher (much higher for $k=2,3$) than in case of {\bf Cube} meshes, even though the global problem is comparable or even smaller in terms of dofs, indicating that the condition number of the stiffness matrix associated to the {\bf Octa} mesh is worse than that of the {\bf Cube} meshes. For $k=1,2$, the Block-Schur preconditioner is completely scalable in terms of iterations, while Block-Reg presents a slight increase. In terms of CPU times, for $k=1$, the most effective solver is Block-Schur, which is between 9 and 30 times as fast as Mumps. For $k=2$ instead, differently from the {\bf Cube} test, Mumps is very effective on {\bf Octa} meshes and results to be the fastest solver, presenting also a good scalability. Note that, with 32 procs and $k=3$, Mumps is 25 and 7 times as fast as the Block-Schur and Block-Reg preconditioners, respectively.

%for both $k=2$ and $k=3$, the performance of the Block-Schur and Block-Reg solvers is comparable with the case of the {\bf CVT} meshes, in terms of GMRES iterations. The solution time are instead significantly smaller, because the global problems have less dofs. Differently from the {\bf CVT} test, Mumps is very effective on {\bf Octa} meshes and results to be the fastest solver, presenting also a good scalability. Note that, with 32 procs and $k=3$, Mumps is 25 and 7 times as fast as the Block-Schur and Block-Reg preconditioners, respectively.

\begin{figure}[!htb]
\begin{center}
\includegraphics[width=0.45\textwidth]{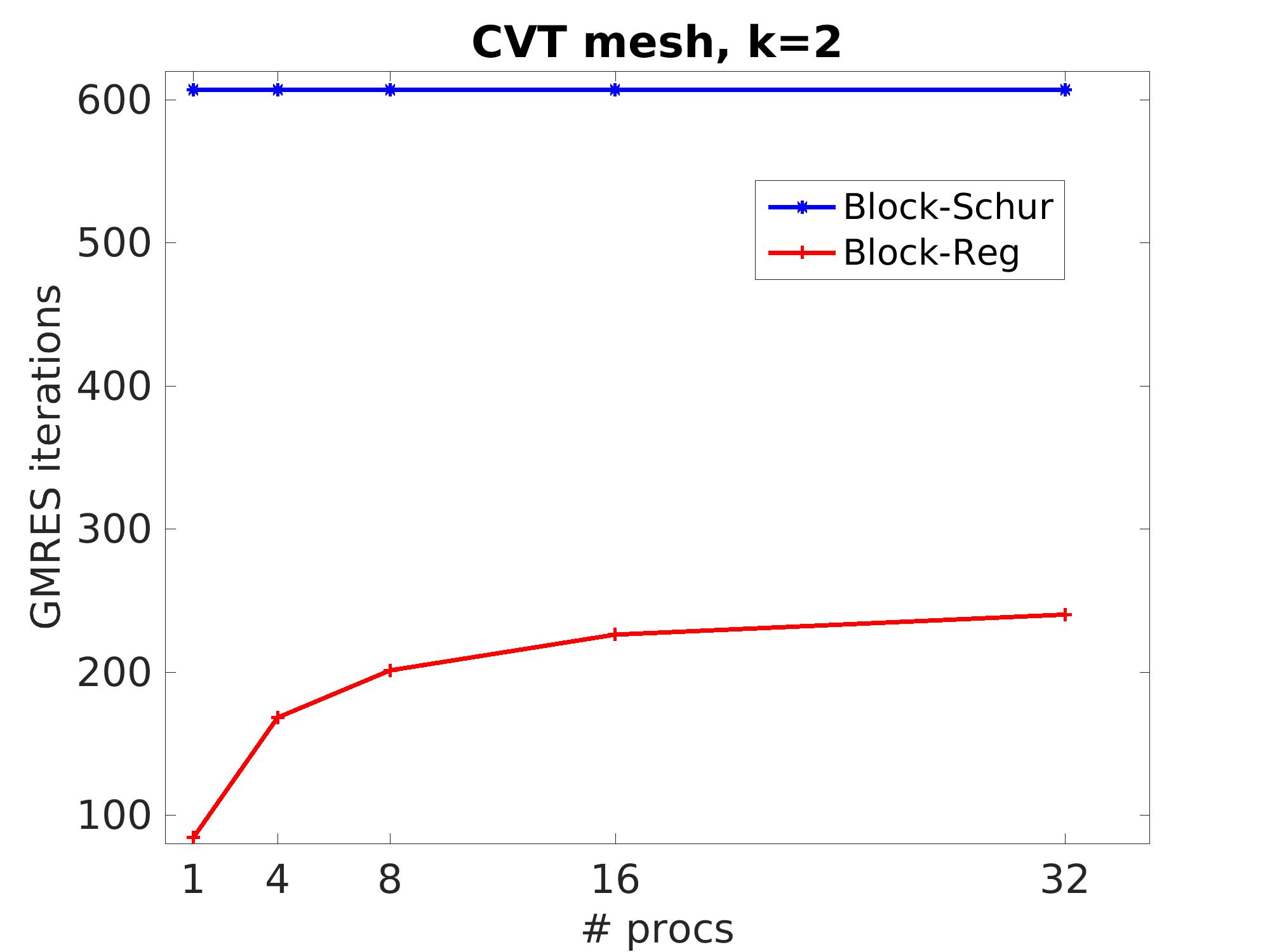}
\includegraphics[width=0.45\textwidth]{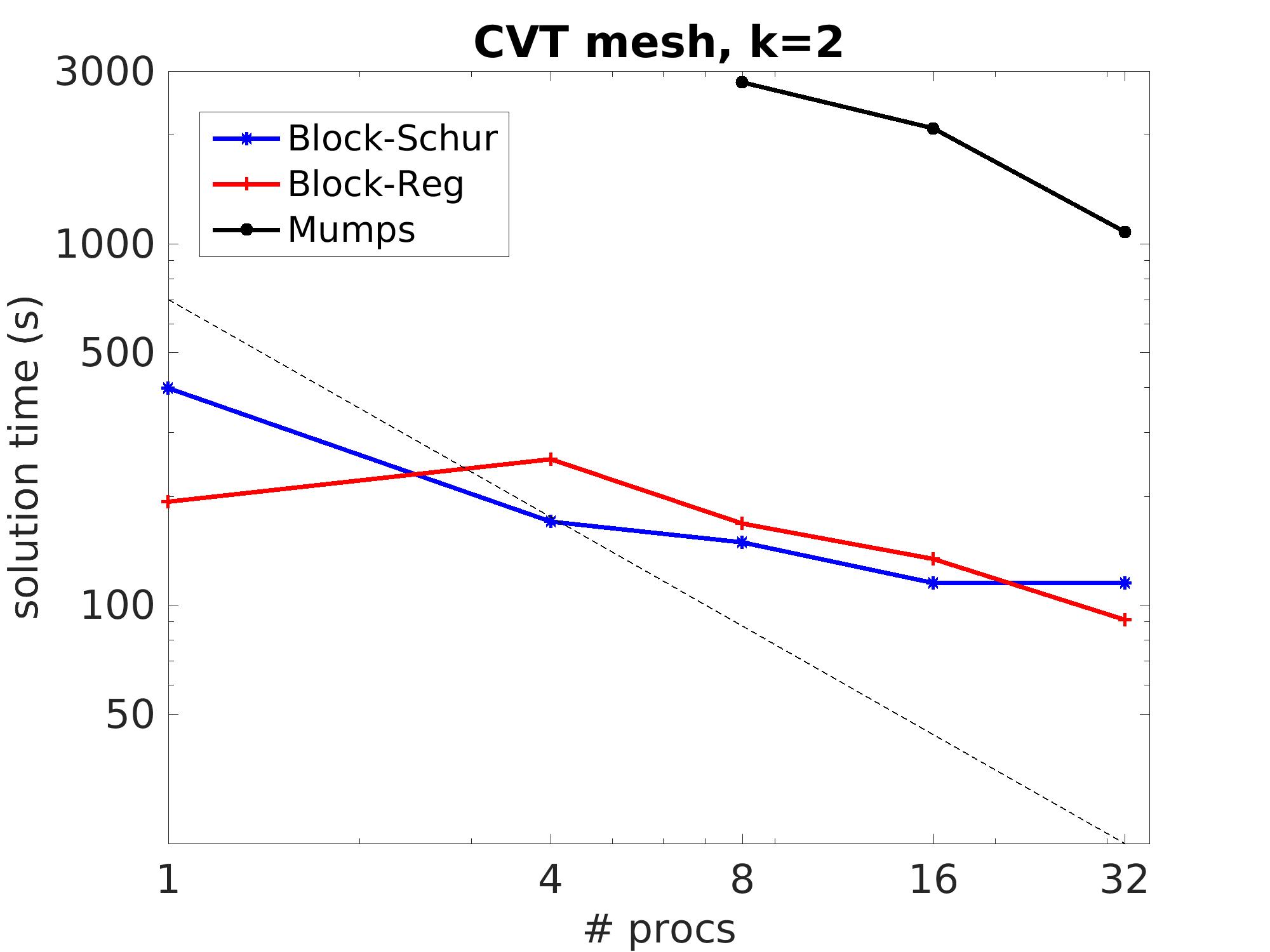}\\
\includegraphics[width=0.45\textwidth]{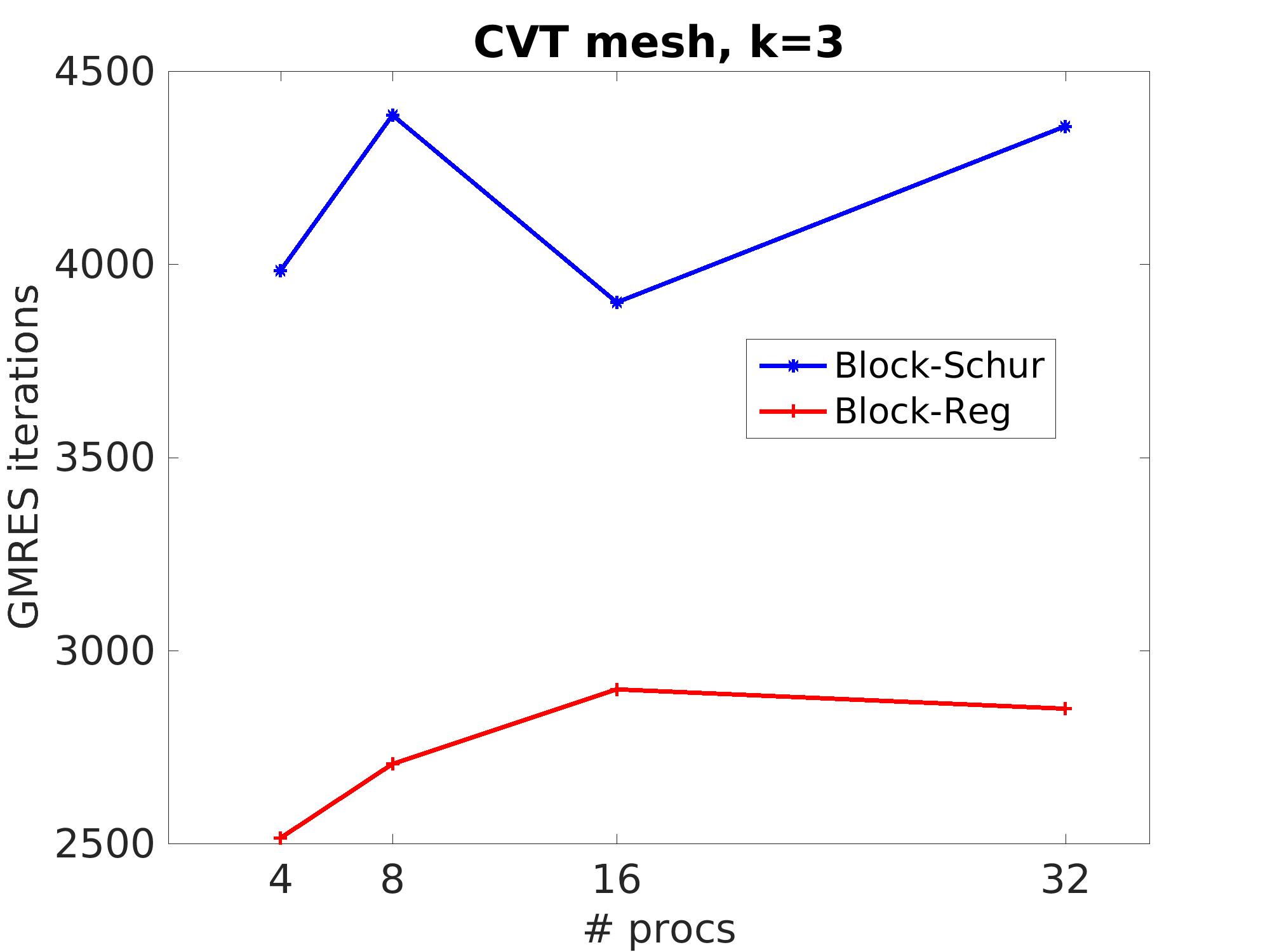}
\includegraphics[width=0.45\textwidth]{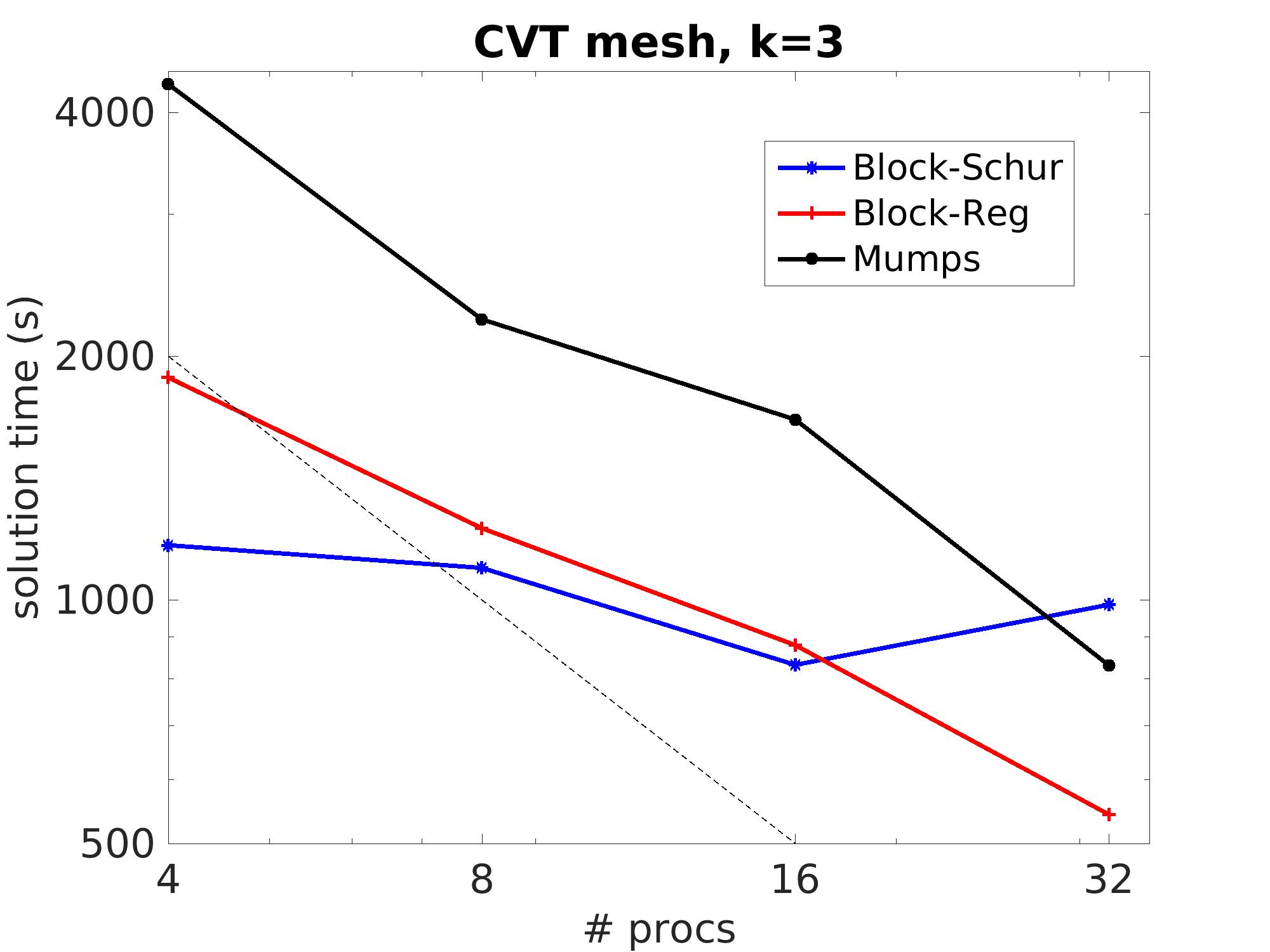}\\
\caption{Strong scaling test, {\bf CVT} meshes. Same format as in Fig. \ref{fig_hexa_scal}.}
\label{fig_voro_scal}
\end{center}
\end{figure}

In case of {\bf CVT} meshes (Table \ref{table_voro_scal} and Fig. \ref{fig_voro_scal}), 
for $k=1,2$ the performance of the Block-Schur and Block-Reg solvers is comparable with the case of the {\bf Octa} meshes, in terms of GMRES iterations.
In terms of CPU time, the speedup values are far from the ideal ones, but the iterative solvers are much more effective than Mumps, which fails for $k=2$ due to Out of Memory for $p=1,\,4$ and is very slow for $p=8,\,16,\,32$. Indeed, for $k=1,2$, the Block-Schur solver results to be between 9 and 300 times as fast as Mumps. For $k=3$ instead, the behavior of the iterative solvers degenerate, due to the severe ill-conditioning of the stiffness matrix. Consequently, Mumps becomes competitive or even faster than the iterative methods.

%In case of {\bf Octa} meshes (Table \ref{table_octa_scal} and Fig. \ref{fig_octa_scal}), for both $k=2$ and $k=3$, the performance of the Block-Schur and Block-Reg solvers is comparable with the case of the {\bf CVT} meshes, in terms of GMRES iterations. The solution time are instead significantly smaller, because the global problems have less dofs. Differently from the {\bf CVT} test, Mumps is very effective on {\bf Octa} meshes and results to be the fastest solver, presenting also a good scalability. Note that, with 32 procs and $k=3$, Mumps is 25 and 7 times as fast as the Block-Schur and Block-Reg preconditioners, respectively.

%&\todo[inline]{Non sarebbe meglio mettere gli esempi in ordine di difficoltà della mesh, ossia 
%{\bf Cube}, {\bf Octa} and {\bf CVT} quindi anche la presentazione dei risultati deve rispecchiare questo ordine. Il pezzo che inizia con ``In case of {\bf Octa} ...'' dovrebbe essere messo prima di ``In case of {\bf CVT}''. Ho visto che hai fatto molti riferimenti nel paragrafo di {\bf Octa} a quelli su {\bf CVT}
%quindi ti chiederei se potessi rivedere un po' questi due paragrafi.}

\begin{table}[!htb]
\begin{center}
\begin{tabular}{c|c|cc|cc|ccc|ccc}
\hline
\multicolumn{12}{c}{{\bf Octa mesh with 30375 elements, k = 1, dofs = 453601}} \\
\hline
$p$ &$S_{p}^{\text{id}}$ &  $T_{ass}$  & $S_{p}$ & \multicolumn{2}{c|}{Mumps} & \multicolumn{3}{c|}{Block-Schur} &  \multicolumn{3}{c}{Block-Reg} \\%[1mm]
&  & & & $T_{sol}$  &$S_{p}$  & it & $T_{sol}$ & $S_{p}$ &it & $T_{sol}$ & $S_{p}$\\
\hline
1       &-	&79	&-	&1901	&-	&79	&60	&-	&93	&61	&-\\
4       &4	&26	&3.0	&551	&3.4	&79	&21	&2.9	&239	&65	&0.9\\
8       &8	&14	&5.6	&336	&5.7	&79	&17	&3.5	&240	&50	&1.2\\
16      &16	&7	&11.3	&205	&9.3	&79	&12	&5.0	&246	&44	&1.4\\
32      &32	&3	&26.3	&116	&16.4	&79	&13	&4.6	&255	&33	&1.8\\
\hline
\hline
\multicolumn{12}{c}{{\bf Octa mesh with 9000 elements, k = 2, dofs = 361201}} \\
\hline
$p$ &$S_{p}^{\text{id}}$ &  $T_{ass}$  & $S_{p}$ & \multicolumn{2}{c|}{Mumps} & \multicolumn{3}{c|}{Block-Schur} &  \multicolumn{3}{c}{Block-Reg} \\%[1mm]
&  & & & $T_{sol}$  &$S_{p}$  & it & $T_{sol}$ & $S_{p}$ &it & $T_{sol}$ & $S_{p}$\\
\hline
1       &-      &342    &-              &326    &-              &559    &276    &-              &113    &109    &-\\
4       &4      &97     &3.5    &107    &3.0    &559    &113    &2.4    &233    &83     &1.3 \\
8       &8      &60     &5.7    &66     &4.9    &559    &102    &2.7    &277    &58     &1.9 \\
16      &16     &32     &10.7   &38     &8.6    &559    &83     &3.3    &274    &61     &1.8 \\
32      &32     &16     &21.4   &24     &13.6   &559    &84     &3.3    &285    &46     &2.4 \\
\hline
\hline
\multicolumn{12}{c}{{\bf Octa mesh with 4608 elements, k = 3, dofs = 378881}} \\
\hline
$p$ &$S_{p}^{\text{id}}$&  $T_{ass}$  & $S_{p}$ & \multicolumn{2}{c|}{Mumps} & \multicolumn{3}{c|}{Block-Schur} &  \multicolumn{3}{c}{Block-Reg} \\%[1mm]
&  & & & $T_{sol}$  &$S_{p}$  & it & $T_{sol}$ & $S_{p}$ &it & $T_{sol}$ & $S_{p}$\\
\hline
4       &-      &683    &-              &108    &-              &2282   &918    &-              &1177   &421    &-\\
8       &2      &354    &1.9    &75     &1.4    &2503   &746    &1.2    &1436   &316    &1.3 \\
16      &4      &192    &3.6    &44     &2.4    &2188   &609    &1.5    &1561   &257    &1.6 \\
32      &8      &98     &7.0    &26     &4.1    &2192   &660    &1.4    &1618   &180    &2.3 \\
\hline
\end{tabular}
\caption{Strong scaling test, Octa meshes. Same format as in Table \ref{table_hexa_scal}.}
\label{table_octa_scal}
\end{center}
\end{table}

\begin{table}[!htb]
\begin{center}
\begin{tabular}{c|c|cc|cc|ccc|ccc}
\hline
\multicolumn{12}{c}{{\bf CVT mesh with 16000 elements, k = 1, dofs = 388043}} \\
\hline
$p$ &$S_{p}^{\text{id}}$&  $T_{ass}$  & $S_{p}$ & \multicolumn{2}{c|}{Mumps} & \multicolumn{3}{c|}{Block-Schur} &  \multicolumn{3}{c}{Block-Reg} \\%[1mm]
&  & & & $T_{sol}$  &$S_{p}$  & it & $T_{sol}$ & $S_{p}$ &it & $T_{sol}$ & $S_{p}$\\
\hline
1       &-	&164	&-	&10687	&-	&138	&35	&-	&72	&80	&-\\
4       &4      &52	&3.1	&2966	&3.6	&138	&15	&2.3	&140	&88	&0.9\\
8       &8      &27	&6.1	&1708	&6.3	&138	&11	&3.2	&162	&63	&1.3\\
16      &16     &14	&11.7	&1230	&8.7	&138	&10	&3.5	&171	&47	&1.7\\
32      &32     &7	&23.4	&614	&17.4	&138	&12	&2.9	&179	&32	&2.5\\
\hline
\hline
\multicolumn{12}{c}{{\bf CVT mesh with 8000 elements, k = 2, dofs = 465721}} \\
\hline
$p$ &$S_{p}^{\text{id}}$&  $T_{ass}$  & $S_{p}$ & \multicolumn{2}{c|}{Mumps} & \multicolumn{3}{c|}{Block-Schur} &  \multicolumn{3}{c}{Block-Reg} \\%[1mm]
&  & & & $T_{sol}$  &$S_{p}$  & it & $T_{sol}$ & $S_{p}$ &it & $T_{sol}$ & $S_{p}$\\
\hline
1       &-      &1688	&-		&F	&-		&607	&398	&-		&84	&193	&-\\
4       &4	&497	&3.4 	&F	&-		&607	&170	&2.3 	&168	&253	&0.8 \\
8       &8	&285	&5.9 	&2793	&-		&607	&149	&2.7 	&201	&168	&1.1 \\
16      &16     &156	&10.8 	&2083	&1.3 	&607	&115	&3.5 	&226	&134	&1.4 \\
32      &32     &83	&20.3 	&1078	&2.6 	&607	&115	&3.5 	&240	&91	&2.1 \\
%64     &64     &43	&39.2 (64)		&899	&3.1 	&607	&396	&1.0 (64)	&706	&517	&1.2 (64)\\
\hline
\hline
\multicolumn{12}{c}{{\bf CVT mesh with 4000 elements, k = 3, dofs = 445951}} \\
\hline
$p$ &$S_{p}^{\text{id}}$ &  $T_{ass}$  & $S_{p}$ & \multicolumn{2}{c|}{Mumps} & \multicolumn{3}{c|}{Block-Schur} &  \multicolumn{3}{c}{Block-Reg} \\%[1mm]
&  & &  & $T_{sol}$  &$S_{p}$  & it & $T_{sol}$ & $S_{p}$ &it & $T_{sol}$ & $S_{p}$\\
\hline
4       &-      &2407	&-		&4335	&-		&3983	&1168	&-		&2514	&1882	&-\\
8       &2      &1276	&1.9 	&2221	&1.9 	&4386	&1095	&1.1 	&2706	&1226	&1.5 \\
16      &4      &779	&3.1 	&1670	&2.6 	&3901	&831	&1.4 	&2899	&879	&2.1 \\
32      &8      &414	&5.8 	&830	&5.2 	&4357	&986	&1.2 	&2849	&543	&3.5 \\
%64     &16 &	&		&	&		&	&	&		&	&	&\\
\hline
\end{tabular}
\caption{Strong scaling test, CVT meshes. Same format as in Table \ref{table_hexa_scal}.}
\label{table_voro_scal}
\end{center}
\end{table}

\subsubsection{Test 2: optimality}

\begin{figure}[!htb]
\begin{center}
\includegraphics[width=0.32\textwidth]{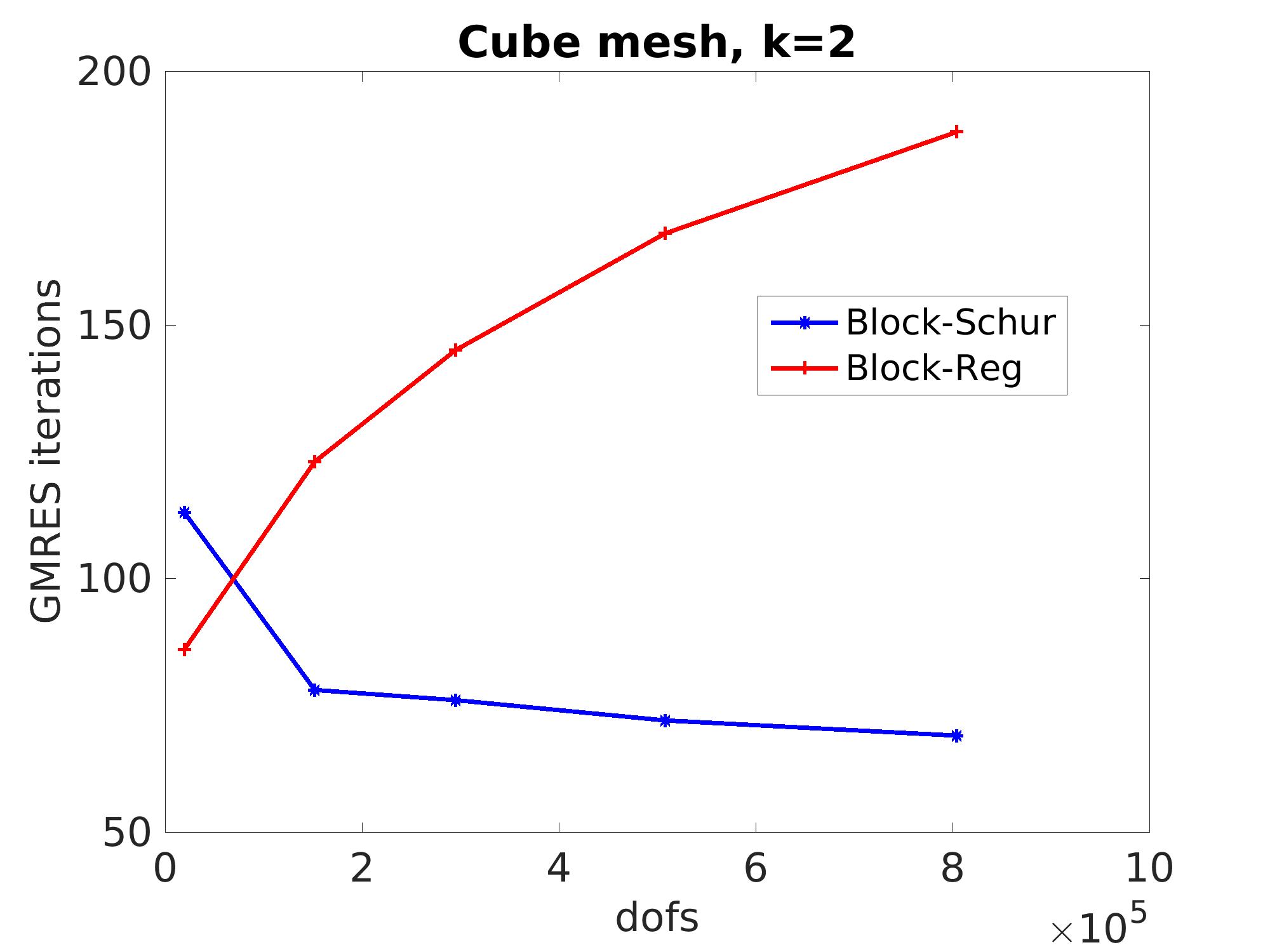}
\includegraphics[width=0.32\textwidth]{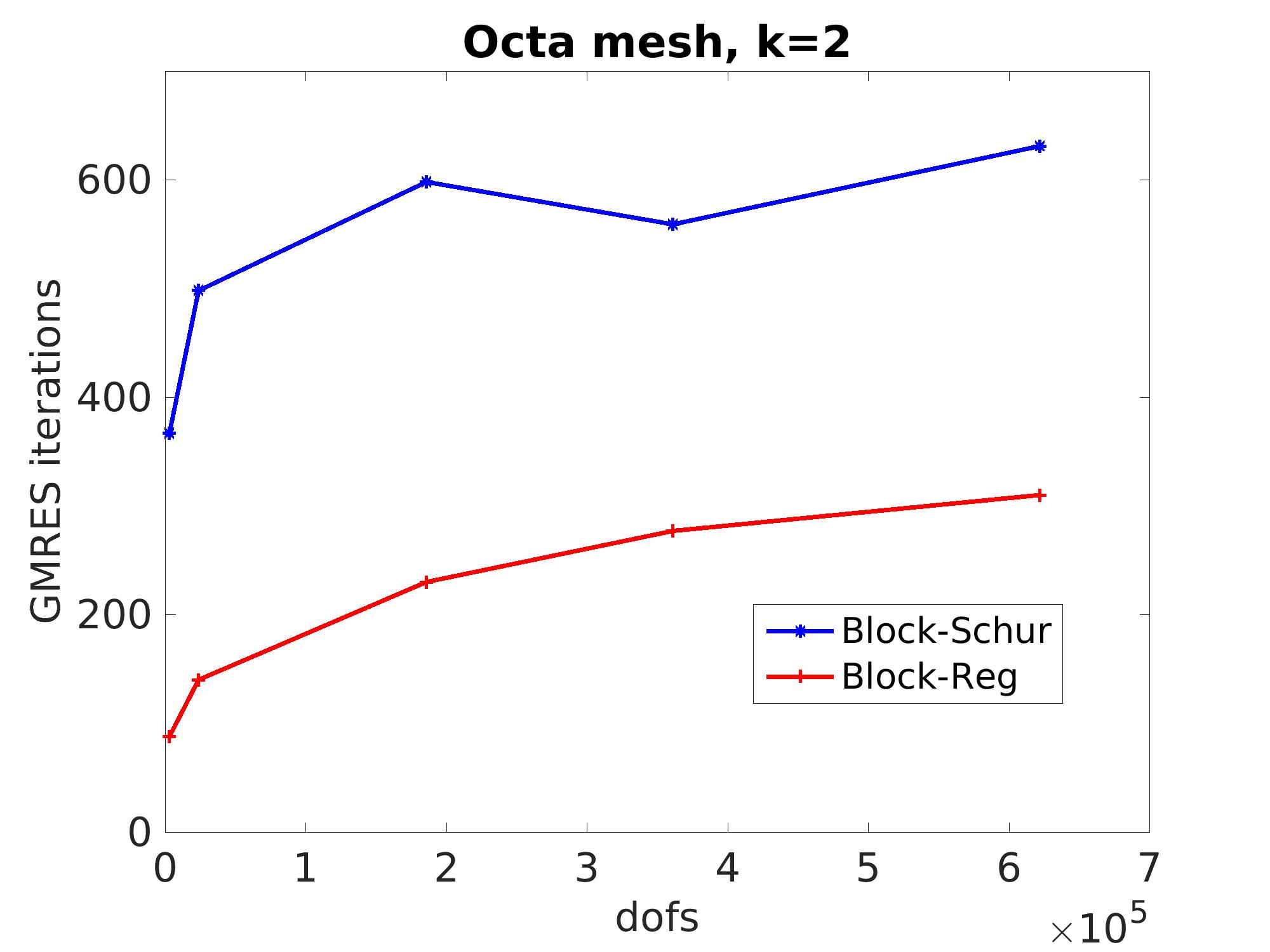}
\includegraphics[width=0.32\textwidth]{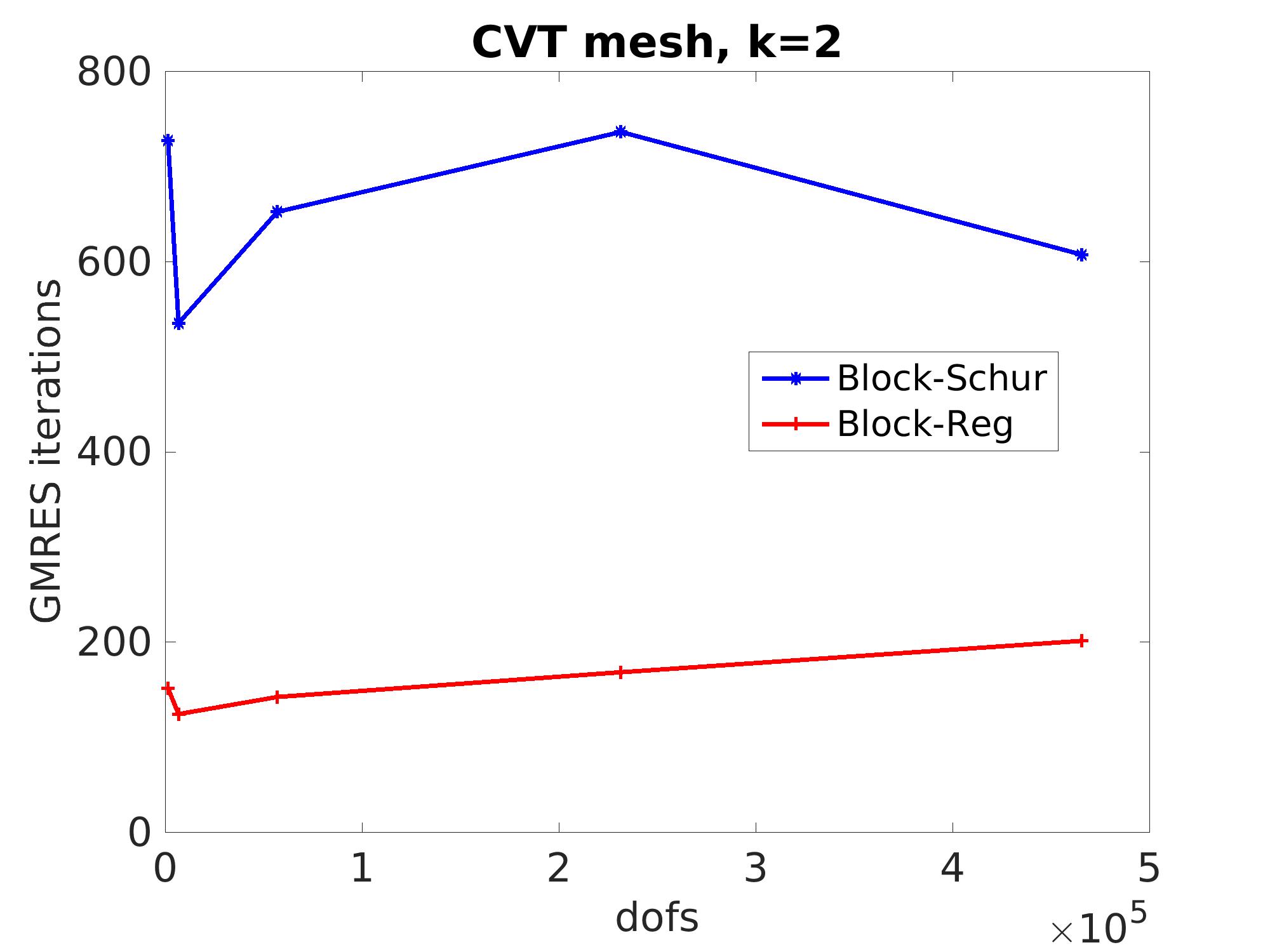}
\caption{Optimality test. GMRES iterations as a function of the number of dofs for $k=2$ VEM discretization on {\bf Cube} (left), {\bf Octa} (center) and {\bf CVT} (right) meshes.}
\label{fig_optim}
\end{center}
\end{figure}

In this test, we investigate the behavior of the parallel solvers when refining the mesh size, thus increasing the number of dofs. The number of processors $p$ and the order of VEM discretization $k$ are kept fixed to 8 and 2, respectively. The results are reported in Tables \ref{table_optim_cube}, \ref{table_optim_voro}, \ref{table_optim_octa} and in Fig. \ref{fig_optim}. In terms of GMRES iterations, the Block-Schur preconditioner exhibits an optimal behavior, since the number of iterations remains bounded from above when varying the number of dofs. The Block-Reg preconditioner instead shows a quasi-optimal behavior, with a growth of GMRES iterations that appears to be logarithmic.

On the largest {\bf Cube} mesh (804385 dofs), the Block-Schur and Block-Reg solvers are about 3 and 7 times as fast as Mumps, respectively. On the largest {\bf Octa} mesh (622081 dofs), Block-Reg results to be the most effective solver. On the largest {\bf CVT} mesh (465721 dofs), the Block-Schur and Block-Reg solvers are about 19 and 17 times as fast as Mumps, respectively.

%\todo[inline]{Manca un commento sui valori delle mesh {\bf CVT}.}

%\todo[inline]{Nelle tabelle di questa sottosezione ho messo la notazione con $\backslash$num per i numeri troppo alti dimmi se ti piace.}

\begin{table}[!htb]
\begin{center}
\begin{tabular}{c|c|c|c|cc|cc}
\hline
\multicolumn{7}{c}{{\bf Cube mesh, k = 2, p=8}} \\
\hline
$N_P$	&dofs	&$T_{ass}$	&Mumps		&\multicolumn{2}{c|}{Block-Schur}	&\multicolumn{2}{c}{Block-Reg}\\%[1mm]
	&  	&	&$T_{sol}$	&it 	&$T_{sol}$			&it 	&$T_{sol}$\\
\hline
\num{512}	&\num{19585}	&1	&1		&113	&1				&86	&1\\
\num{4096}	&\num{152065}	&13	&14		&78	&9				&123	&9\\
\num{8000}	&\num{295201}	&29	&73		&76	&29				&145	&22\\	
\num{13824}	&\num{508033}	&53	&270		&72	&77				&168	&48\\
\num{21952}	&\num{804385}	&84	&604		&69	&217				&188	&87\\
\hline
\end{tabular}
\caption{Optimality test, {\bf Cube} meshes. $N_P$ = number of mesh polyhedra; dofs = degrees of freedom; $T_{ass}$:=assembling time in seconds; $T_{sol}$:=solution time in seconds; it:=GMRES iterations. The order of VEM discretization is fixed to k=2. The runs are performed on 8 processors (p=8).}
\label{table_optim_cube}
\end{center}
\end{table}

\begin{table}[!htb]
\begin{center}
\begin{tabular}{c|c|c|c|cc|cc}
\hline
\multicolumn{7}{c}{{\bf Octa mesh, k = 2, p=8}} \\
\hline
$N_P$     &dofs   &$T_{ass}$      &Mumps          &\multicolumn{2}{c|}{Block-Schur}       &\multicolumn{2}{c}{Block-Reg}\\%[1mm]
        &       &       &$T_{sol}$      &it     &$T_{sol}$                      &it     &$T_{sol}$\\
\hline
\num{72}        &\num{3121}	&$<1$	&$<1$		&367	&2				&88	&$<1$\\
\num{576}	&\num{23809}	&3	&$<1$		&498	&4				&140	&2\\
\num{4608}	&\num{185857}	&25	&13		&598	&35				&230	&24\\
\num{9000}	&\num{361201}	&60	&66		&559	&102				&277	&58\\
\num{15552}	&\num{622081}	&111	&190		&631	&283				&310	&113\\
\hline
\end{tabular}
\caption{Optimality test, {\bf Octa} meshes. Same format as in Table \ref{table_optim_cube}.}
\label{table_optim_octa}
\end{center}
\end{table}

\begin{table}[!htb]
\begin{center}
\begin{tabular}{c|c|c|c|cc|cc}
\hline
\multicolumn{7}{c}{{\bf CVT mesh, k = 2, p=8}} \\
\hline
$N_P$     &dofs   &$T_{ass}$      &Mumps          &\multicolumn{2}{c|}{Block-Schur}       &\multicolumn{2}{c}{Block-Reg}\\%[1mm]
        &       &       &$T_{sol}$      &it     &$T_{sol}$                      &it     &$T_{sol}$\\
\hline
\num{27}        &\num{1495}     &1      &$<1$           &727    &2                              &151    &$<1$\\
\num{125}       &\num{6841}     &3      &$<1$           &535    &2                              &124    &2\\
\num{1000}      &\num{57019}    &33     &20             &652    &10                             &142    &16\\
\num{4000}      &\num{231571}   &143    &467            &736    &46                             &168    &75\\
\num{8000}      &\num{465721}   &285    &2793           &607    &149                            &201    &168\\
\hline
\end{tabular}
\caption{Optimality test, {\bf CVT} meshes. Same format as in Table \ref{table_optim_cube}.}
\label{table_optim_voro}
\end{center}
\end{table}

%-------------------------------------------------------------

\section{Conclusion}

In this paper we constructed and numerically analyzed a general order VEM approximation scheme for three-dimensional scalar elliptic equations in mixed form. 
The convergence tests have demonstrated the effectiveness and the robustness of the proposed algorithms on different polyhedral grids.
Moreover, we developed a parallel solver for the solution of the saddle point linear systems arising from the discretization process, by exploiting both direct and iterative parallel solution methods, preconditioned by block-diagonal preconditioners. 
The numerical tests performed on a Linux cluster have shown that the proposed iterative methods are more effective than the Mumps direct solver for low order ($k=1,2$) discretizations, while, using high order discretizations, they suffer due to the severe ill-conditioning of the linear system matrix and Mumps results to be the fastest solver. 
%\todo[inline]{Mi sembra che questa frase messa così sminuisca il nostro lavoro...
%Possiamo metterla più in positivo?}
More research is needed in future to contruct preconditioners resulting to be robust with respect to high order of VEM discretizations. Further developments of the present investigation might be the construction and analysis of Balancing Domain Decomposition by Constraints (BDDC) preconditioners for VEM discretizations of elliptic equations in mixed form and the extension of the parallel solvers studied here to the solution of other saddle point linear systems arising from VEM discretizations of Stokes or Maxwell equations.

%-------------------------------------------------------------

\section*{Acknowledgments}

The authors would like to acknowledge INDAM-GNCS for the support.
Moreover they would like to thank Lourenco Beir{\~a}o~da Veiga and Alessandro Russo for 
many helpful discussions and suggestions.

%-------------------------------------------------------------
% References 
%-------------------------------------------------------------
\bibliographystyle{plain}
\bibliography{VEM}

\end{document}